\documentclass[a4paper,11pt]{article}
\usepackage{amsmath}
\usepackage{amssymb}
\usepackage{amsthm}
\usepackage{bm}
\usepackage[margin = 3cm]{geometry}
\usepackage{hyperref}
\usepackage{enumitem}

\newtheorem{theorem}{Theorem}[section]

\newtheorem{assump}{Assumption}

\theoremstyle{definition}
\newtheorem{definition}[theorem]{Definition}
\newtheorem{remark}{Remark}

\numberwithin{equation}{section}

\begin{document}

\title{Analysis of a Cahn--Hilliard system with non-zero Dirichlet conditions modeling tumor growth with chemotaxis}

\author{Harald Garcke \footnotemark[1] \and Kei Fong Lam \footnotemark[1]}

\date{\today }

\maketitle

\renewcommand{\thefootnote}{\fnsymbol{footnote}}
\footnotetext[1]{Fakult\"at f\"ur Mathematik, Universit\"at Regensburg, 93040 Regensburg, Germany
({\tt \{Harald.Garcke, Kei-Fong.Lam\}@mathematik.uni-regensburg.de}).}

\begin{abstract}
We consider a diffuse interface model for tumor growth consisting of a Cahn--Hilliard equation with source terms coupled to a reaction-diffusion equation, which models a tumor growing in the presence of a nutrient species and surrounded by healthy tissue.  The well-posedness of the system equipped with Neumann boundary conditions was found to require regular potentials with quadratic growth.  In this work, Dirichlet boundary conditions are considered, and we establish the well-posedness of the system for regular potentials with higher polynomial growth and also for singular potentials.  New difficulties are encountered due to the higher polynomial growth, but for regular potentials, we retain the continuous dependence on initial and boundary data for the chemical potential and for the order parameter in strong norms as established in the previous work.  Furthermore, we deduce the well-posedness of a variant of the model with quasi-static nutrient by rigorously passing to the limit where the ratio of the nutrient diffusion time-scale to the tumor doubling time-scale is small.
\end{abstract}

\noindent \textbf{Key words. } Tumor growth; phase field model; Cahn--Hilliard equation; reaction-diffusion equations; chemotaxis; weak solutions; Dirichlet boundary conditions; well-posedness; asymptotic analysis; singular potentials. \\

\noindent \textbf{AMS subject classification. } 35D30, 35Q92, 35K57, 35B40, 92C17.

\section{Introduction}
We consider the following system of equations describing a two component mixture of tumor tissue and healthy (surrounding) tissue in the presence of a chemical species acting as nutrient for the tumor,
\begin{subequations}\label{Intro:CHN}
\begin{alignat}{2}
\partial_{t}\varphi & = \Delta \mu +  (\lambda_{p} \sigma - \lambda_{a})h(\varphi) && \text{ in } \Omega \times (0,T), \label{Intro:varphi} \\
\mu & = \frac{\gamma}{\varepsilon} \Psi'(\varphi) - \gamma \varepsilon \Delta \varphi - \chi \sigma && \text{ in } \Omega \times (0,T), \label{Intro:mu} \\
\kappa \partial_{t} \sigma & = \mathrm{div} \, ( D(\varphi) (\nabla \sigma - \eta \nabla \varphi)) - \lambda_{c}  \sigma h(\varphi) && \text{ in } \Omega \times (0,T), \label{Intro:sigma} \\
\sigma & = \sigma_{\infty}, \quad \varphi = -1, \quad \mu = \mu_{\infty} && \text{ on } \Gamma \times (0,T), \\
\varphi(x,0)& = \varphi_{0}(x), \quad \sigma(x,0) = \sigma_{0}(x) && \text{ in } \Omega.
\end{alignat}
\end{subequations}
Here, $T > 0$ denotes a fixed time and $\Omega \subset \mathbb{R}^{3}$ denotes a bounded domain with boundary $\Gamma := \partial \Omega$.  The system \eqref{Intro:CHN} describes a diffuse interface model for tumor growth via a Cahn--Hilliard equation coupled with a reaction-diffusion equation for a nutrient species, whose concentration we denote as $\sigma$.  The order parameter which distinguishes the two components is denoted by $\varphi$, where the region $\{\varphi = 1\}$ corresponds to the tumor phase and the region  $\{\varphi = -1\}$ corresponds to the healthy tissue phase.

The function $\Psi'$ is the derivative of a potential $\Psi$ that has two equal minima at $\pm 1$, $\varepsilon > 0$ is a parameter related to the interfacial thickness, $\gamma > 0$ denotes the surface tension, $\lambda_{p}$, $\lambda_{a}$, and $\lambda_{c}$ are non-negative constants denoting the tumor proliferation rate, tumor apoptosis rate, and the nutrient consumption rate, respectively.  The positive mobility $D(\varphi)$ corresponds to the diffusivity of the nutrient, and $h(\varphi)$ is an interpolation function such that $h(-1) = 0$ and $h(1) = 1$, with the simplest example being $h(\varphi) = \frac{1}{2}(1+\varphi)$.  Here, $\kappa \geq 0$ is a constant such that, for a constant mobility $D(\varphi) = D_{0}$, the ratio $\kappa/D_{0}$ represents the ratio of the nutrient diffusion time-scale and the tumor doubling time-scale.  Finally, $\chi$ and $\eta$ are non-negative constants representing chemotaxis (movement of tumor cells towards regions of high nutrients) and active transport (establishment of persistent nutrient concentration gradient in the vicinity of the tumor interface), respectively.  We refer the reader to \cite{GLSS} for more details regarding the derivation of the model and the modeling of the chemotaxis and active transport mechanisms.

We point out that the well-posedness of \eqref{Intro:CHN} with $\kappa = 1$, $h(\cdot) \in C^{0}(\mathbb{R}) \cap L^{\infty}(\mathbb{R})$, homogeneous Neumann boundary conditions for $\mu$ and $\varphi$, and Robin boundary condition for $\sigma$ has been studied by the authors in \cite{GLNeu}.  In contrast, here we consider Dirichlet boundary conditions for the following reasons.
\begin{itemize}
\item In \cite{GLNeu}, the class of admissible potential $\Psi$ is restricted to potentials with at most quadratic growth.  Here, for the existence of weak solutions, we allow $\Psi$ to have arbitrary polynomial growth (see \eqref{assump:Psi} below), and for regularity and continuous dependence on initial data, we are restricted to the class of potentials that have polynomial growth of order up to (but not including) 6 in dimension $d = 3$ (see \eqref{assump:Psi''} and \eqref{Psi'':uniqueness}).
\item In \cite{GLNeu}, we are unable to pass to the limit $\kappa \to 0$ to deduce if weak solutions of \eqref{Intro:CHN} converge (in some appropriate sense) to weak solutions of the following quasi-static model
\begin{subequations}\label{Intro:quasistatic}
\begin{alignat}{3}
\partial_{t} \varphi & = \Delta \mu + (\lambda_{p} \sigma - \lambda_{a}) h(\varphi) && \text{ in } \Omega \times (0,T),\\
\mu & = \frac{\gamma}{\varepsilon} \Psi'(\varphi) - \gamma \varepsilon \Delta \varphi - \chi \sigma && \text{ in } \Omega \times (0,T), \\
0 & = \mathrm{div} \, ( D(\varphi) (\nabla \sigma - \eta \nabla \varphi)) - \lambda_{c} \sigma h(\varphi)  && \text{ in } \Omega \times (0,T), \\
\sigma & = \sigma_{\infty}, \quad \varphi = -1, \quad \mu = \mu_{\infty} && \text{ on } \Gamma \times (0,T), \\
\varphi(x,0)& = \varphi_{0}(x) && \text{ in } \Omega.
\end{alignat}
\end{subequations}
The well-posedness of \eqref{Intro:quasistatic} is proved separately from \eqref{Intro:CHN} in \cite{GLNeu}.  However, with Dirichlet boundary conditions, we prove that the quasi-static model can be realized as the limit system from \eqref{Intro:CHN} as $\kappa \to 0$ (see Theorem \ref{thm:truncated:quasi} below).
\end{itemize}

We briefly discuss the motivation for the Dirichlet boundary conditions.  For simplicity, in taking $\kappa = \gamma = \varepsilon = \lambda_{p} = \lambda_{c} = D(\varphi)= 1$, $\lambda_{a} = 0$, $\eta = \chi$, and zero Neumann boundary conditions in \eqref{Intro:CHN}, there is a natural energy identity of the form
\begin{equation}\label{Neumann:Energy}
\begin{aligned}
& \frac{d}{dt} \int_{\Omega} \left ( \Psi(\varphi) + \frac{1}{2}|\nabla \varphi|^{2} + \frac{1}{2} |\sigma|^{2} + \chi \sigma (1-\varphi) \right ) \, dx \\
& \qquad + \int_{\Omega} |\nabla \mu|^{2} + |\nabla (\sigma - \chi \varphi)|^{2} \, dx 
\end{aligned}
\end{equation}\begin{equation*}
\begin{aligned}& \quad = \int_{\Omega} h(\varphi) \sigma \mu - h(\varphi) \sigma (\sigma + \chi(1-\varphi)) \, dx.
\end{aligned}
\end{equation*}
The main difficulty in deriving the first useful a priori estimate from \eqref{Neumann:Energy} lies in controlling the product $h(\varphi) \sigma \mu$ on the right-hand side with the left-hand side.  In the absence of any a priori knowledge about the integrability of $\sigma$, it seems natural to use the Cauchy--Schwarz inequality, which leads to the appearance of a term $\|\mu \|_{L^{2}}^{2}$ on the right-hand side.  To control this with $\|\nabla \mu \|_{L^{2}}^{2}$ on the left-hand side with the Poincar\'{e} inequality, we have to obtain an estimate for the square of the mean of $\mu$, which is related to the square of the mean of $\Psi'(\varphi)$, leading to a integral inequality of the form
\begin{align}\label{quadratic:growth}
 \int_{\Omega} \left ( \Psi(\varphi) + |\nabla \varphi|^{2} + |\sigma|^{2} \right ) \, dx \leq C \int_{0}^{s} \left (1 + \| \Psi'(\varphi) \|_{L^{1}}^{2} + \| \nabla \varphi \|_{L^{2}}^{2} + \| \sigma \|_{L^{2}}^{2} \right ) \, dt.
\end{align}

This is in contrast with the standard Cahn--Hilliard equation, where the first a priori estimate can be derived without the need to estimate the mean of $\mu$.  A scaling argument for the polynomial growth of $\Psi$ (cf. \cite[\S 7]{GLNeu}) shows that in order to satisfy \eqref{quadratic:growth} the potential $\Psi$ can have at most quadratic growth, which does not cover the classical quartic potential $W(s) := (s^{2}-1)^{2}$ that commonly appears in phase field models.  Therefore, our current consideration of Dirichlet boundary condition for $\mu$ seeks to elevate this restriction on the growth of $\Psi$, which bypasses the need to estimate the mean of $\mu$ with the use of the Poincar\'{e} inequality.

Furthermore, to rigorously obtain the quasi-static limit $\kappa \to 0$ in \eqref{Intro:CHN}, we have to control the source terms on the right-hand side of \eqref{Neumann:Energy} involving $\sigma$ without making use of $\frac{d}{dt} \kappa \|\sigma\|_{L^{2}}^{2}$ that will appear on the left-hand side.  This can be achieved again through the use of the Poincar\'{e} inequality and Dirichlet boundary condition for $\sigma$.

In Section \ref{sec:aprioriest} we will derive the following energy integral identity for \eqref{Intro:CHN}:  Let $\omega = \sigma - \sigma_{\infty}$ and $\mathcal{X} > 0$ be a positive constant yet to be determined, then we have for $s \in (0,T]$,
\begin{align*}
& \int_{\Omega} \left [ \frac{\gamma}{\varepsilon} \Psi(\varphi) + \frac{\gamma \varepsilon}{2} |\nabla \varphi|^{2} + \mathcal{X} \frac{\kappa}{2} |\omega|^{2} \right ](s) \, dx \\
& \qquad + \int_{0}^{s} \int_{\Omega} |\nabla \mu|^{2} + \mathcal{X} D(\varphi) |\nabla \omega|^{2} + \lambda_{c}h(\varphi) |\omega|^{2} \, dx \, dt \\
& \quad = \int_{\Omega} \left [ \frac{\gamma}{\varepsilon} \Psi(\varphi) + \frac{\gamma \varepsilon}{2} |\nabla \varphi|^{2} + \mathcal{X} \frac{\kappa}{2} |\omega|^{2} \right ](0) \,dx \\
& \qquad + \int_{0}^{s} \int_{\Omega} \partial_{t}\varphi (\mu_{\infty} + \chi \sigma_{\infty}) - \chi \nabla \mu \cdot \nabla \omega + \nabla \mu \cdot \nabla \mu_{\infty} \, dx \, dt \\
& \qquad + \int_{0}^{s} \int_{\Omega} h(\varphi)(\lambda_{p} \sigma - \lambda_{a})(\mu - \mu_{\infty} + \chi \omega) \, dx \, dt \\
& \qquad + \mathcal{X} \int_{0}^{s} \int_{\Omega} D(\varphi) \eta \nabla \varphi \cdot \nabla \omega - \kappa \partial_{t} \sigma_{\infty} \omega - D(\varphi) \nabla \sigma_{\infty} \cdot \nabla \omega - \lambda_{c} h(\varphi) \sigma_{\infty} \omega \, dx \, dt.
\end{align*}

This is obtained from testing \eqref{Intro:varphi} with $\mu - \mu_{\infty} + \chi (\sigma - \sigma_{\infty})$, \eqref{Intro:mu} with $\partial_{t}\varphi$ and \eqref{Intro:sigma} with $\mathcal{X} (\sigma - \sigma_{\infty})$.  The boundary conditions $\mu_{\infty}$ and $\sigma_{\infty}$ appear due to the fact that we can only test with functions that have zero trace on $\Gamma$ in the weak formulation of \eqref{Intro:CHN}, and from the above energy identity, assumptions on the spatial gradients and time derivatives (in order to handle the term $\partial_{t}\varphi (\mu_{\infty} + \chi \sigma_{\infty})$) for $\mu_{\infty}$ and $\sigma_{\infty}$, such as \eqref{assump:boundarydata} below, will be needed to derive useful a priori estimates.  Similar assumptions for the boundary values have also been used in the seminal work of Alt and Luckhaus \cite{AltLuck}.

Let us mention that, while a Dirichlet boundary condition for $\sigma$ is reasonable from the modeling viewpoint, which represents a source of nutrients into the domain $\Omega$, there has been little consideration of a Dirichlet boundary condition for the chemical potential $\mu$ in the literature.  To the authors' best knowledge, the Dirichlet boundary condition $\mu = 0$ has been considered in \cite{Chen, Frieboes, Wise}, and this has the interpretation that the cells are allowed to flow freely across the boundary $\Gamma$.  In this work, we treat a more general boundary value $\mu_{\infty}$ and investigate the assumptions that are sufficient for well-posedness.

It is also possible to consider more general boundary conditions $\varphi_{\infty}$ for $\varphi$, but in this work we restrict to the case $\varphi_{\infty} = -1$ which allows for the physical interpretation that the tumor region is enclosed by the healthy tissue.  Using the Yosida approximation and due to a priori estimates which do not depend on the Yosida parameter, we can extend our analysis to potentials of the form
\begin{align}\label{defn:singularpot}
\Psi(y) = \hat{\beta}(y) + \Lambda(y),
\end{align}
where $\hat{\beta}: \mathbb{R} \to [0,\infty]$ is a convex, proper (i.e., not identically $\infty$), lower semicontinuous function and $\Lambda : \mathbb{R} \to \mathbb{R}$ is a $C^{1,1}$ perturbation with at most quadratic growth.  It is possible that $\hat{\beta}$ does not possess a classical derivative, and so $\Psi'$ may not be well-defined.  But we can consider the subdifferential $\partial \hat{\beta}$ as a notion of generalized derivative.

We briefly recall that the effective domain of $\hat{\beta}$ is defined as $D(\hat{\beta}) := \{ x \in \mathbb{R} : \hat{\beta}(x) < \infty \}$, and the effective domain for a possibly multivalued mapping $T:\mathbb{R} \to 2^{\mathbb{R}}$ is $D(T) := \{ x \in \mathbb{R} : Tx \neq \emptyset \}$.  The subdifferential $\partial \hat{\beta} :\mathbb{R} \to 2^{\mathbb{R}}$ of $\hat{\beta}$ is a potentially multivalued mapping defined as
\begin{align}\label{subdifferential}
\partial \hat{\beta}(x) := \{ f \in \mathbb{R} : \hat{\beta}(y) - \hat{\beta}(x) \geq f(y-x) \; \forall y \in \mathbb{R} \} \text{ for } x \in \mathbb{R},
\end{align}
with domain $D(\partial \hat{\beta}) = \{ x \in \mathbb{R} : \partial \hat{\beta}(x) \neq \emptyset \}$.  Using the notation $\beta(x) := \partial \hat{\beta}(x)$, we denote by $\beta^{0}(x)$ the element of the set $\beta(x)$ such that
\begin{align*}
|\beta^{0}(x)| :=  \inf_{f \in \beta(x)} |f| \text{ for } x \in D(\beta).
\end{align*}

The well-posedness theory of \eqref{Intro:CHN} will allow us to deduce the existence and uniqueness of weak solutions to
\begin{subequations}\label{singular:CHN}
\begin{alignat}{2}
\partial_{t}\varphi & = \Delta \mu + (\lambda_{p} \sigma - \lambda_{a} ) h(\varphi) && \text{ in } \Omega \times (0,T), \\
\mu & = \frac{\gamma}{\varepsilon} (\psi + \Lambda'(\varphi)) - \gamma \varepsilon \Delta \varphi - \chi \sigma && \text{ in } \Omega \times (0,T), \\
\kappa \partial_{t} \sigma & = \mathrm{div} \, (D(\varphi) (\nabla \sigma - \eta \nabla \varphi)) - \lambda_{c} \sigma h(\varphi) && \text{ in } \Omega \times (0,T), \\
\sigma & = \sigma_{\infty}, \quad \varphi = -1, \quad \mu = \mu_{\infty} && \text{ on } \Gamma \times (0,T), \\
\varphi(x,0)& = \varphi_{0}(x), \quad \sigma(x,0) = \sigma_{0}(x) && \text{ in } \Omega,
\end{alignat}
\end{subequations}
where $\psi \in \beta(\varphi) = \partial \hat{\beta}(\varphi)$ denotes a selection, as $\beta(\varphi)$ can be multivalued.  In particular, if we consider $\Psi$ to be a potential of double-obstacle type
\begin{align}\label{DoubleObstaclebeta}
\hat{\beta}(y) = \begin{cases}
0 & \text{ for } y \in [-1,1],\\
\infty & \text{ otherwise},
\end{cases} \quad \beta(y) = \partial \hat{\beta}(y) = \begin{cases}
(-\infty,0] & \text{ if } y = -1, \\
0 & \text{ if } y \in (-1,1), \\
[0,\infty) & \text{ if } y = 1,
\end{cases}
\end{align}
and
\begin{align}\label{DoubleObstacleLambda}
 \Lambda(y) = \frac{1}{2} (1-y^{2}),
\end{align}
then $\varphi \in [-1,1] = D(\beta)$ a.e. in $\Omega \times (0,T)$.

Lastly, we compare with some of the recent well-posedness results on phase field models for tumor growth.  The model studied in \cite[(1.1) - (1.4)]{CGH} (with $\alpha = 0$) and \cite[(1.1) - (1.4)]{FGR} bears the most resemblance to \eqref{Intro:CHN},
\begin{subequations}\label{ModelRocca}
\begin{alignat}{3}
\partial_{t}\varphi &= \Delta \mu + h(\varphi)(\sigma - \mu) && \text{ in } \Omega \times (0,T), \label{Model:Rocca:varphi} \\
\mu & = \Psi'(\varphi) - \Delta \varphi && \text{ in } \Omega \times (0,T), \label{Model:Rocca:mu} \\
\partial_{t}\sigma & = \Delta \sigma - h(\varphi)(\sigma - \mu) && \text{ in } \Omega \times (0,T), \label{Model:Rocca:sigma} \\
0 & = \nabla \varphi \cdot \nu = \nabla \mu \cdot \nu = \nabla \sigma \cdot \nu && \text{ on } \Gamma \times (0,T).
\end{alignat}
\end{subequations}
The source term $h(\varphi)(\sigma - \mu)$ appearing in \eqref{Model:Rocca:varphi} and \eqref{Model:Rocca:sigma} is motivated by linear phenomenological constitutive laws for chemical reactions (see \cite{Hawkins} for more details), and is different compared to our choice of source terms in \eqref{Intro:CHN}.  The well-posedness of the system \eqref{ModelRocca} has been established in \cite{CGH,FGR} for large classes of $\Psi$ and $h$.

Another class of models that describes tumor growth uses a Cahn--Hilliard--Darcy system
\begin{alignat*}{3}
\mathrm{div} \, \boldsymbol{v} & = \mathcal{S} && \text{ in } \Omega \times (0,T), \\
\boldsymbol{v} &= -\nabla p + \mu \nabla \varphi && \text{ in } \Omega \times (0,T), \\
\partial_{t} \varphi + \mathrm{div} \, (\varphi \boldsymbol{v}) & = \Delta \mu + \mathcal{S} && \text{ in } \Omega \times (0,T), \\
\mu & = \Psi'(\varphi) - \Delta \varphi && \text{ in } \Omega \times (0,T), \\
0 & = \nabla \varphi \cdot \nu = \nabla \mu \cdot \nu = \boldsymbol{v} \cdot \nu && \text{ on } \Gamma \times (0,T),
\end{alignat*}
where $\boldsymbol{v}$ denote a mixture velocity, $p$ denotes the pressure, and $\mathcal{S}$ is a mass exchange term.  The existence of strong solutions in 2D and 3D have been studied in \cite{LTZ} for the case $\mathcal{S} = 0$.  For the case where $\mathcal{S} \neq 0$ is prescribed, global weak solutions and unique local strong solutions in 3D, as well as global strong well-posedness and long-time behavior in 2D  can be found in \cite{JWZ}.  We also mention the work of \cite{Bosia} which treats a related system known as the Cahn--Hilliard--Brinkman system.  For the Cahn--Hilliard--Darcy system with nutrient which was proposed in \cite{GLSS}, the global existence of weak solutions have been established in \cite{GLDarcy}, and also in \cite{GLRome} for the quasi-static variant.

The plan for the rest of this paper is as follows.  In Section \ref{sec:mainresults}, we state the assumptions and the results for \eqref{Intro:CHN}, \eqref{Intro:quasistatic} and \eqref{singular:CHN}.  In Section \ref{sec:regularPot} we perform a Galerkin procedure to deduce existence of weak solutions to \eqref{Intro:CHN}, and show further regularity and continuous dependence on initial and boundary data.  In Section \ref{sec:quasistaticlimit} we pass to the limit $\kappa \to 0$ in \eqref{Intro:CHN} and deduce the existence result for \eqref{Intro:quasistatic}, and then we prove further regularity properties and continuous dependence on initial and boundary data.  The regular weak solutions to \eqref{Intro:CHN} satisfy an energy inequality and this will allow us to employ the Yosida approximation to prove the existence of weak solutions to \eqref{singular:CHN}, which is done in Section \ref{sec:singularPot}.

\subsection*{Notation} For convenience, we use the notation $L^{p} := L^{p}(\Omega)$ and $W^{k,p} := W^{k,p}(\Omega)$ for any $p \in [1,\infty]$, $k > 0$ to denote the standard Lebesgue spaces and Sobolev spaces equipped with the norms $\| \cdot \|_{L^{p}}$ and $\| \cdot \|_{W^{k,p}}$.  In the case $p = 2$ we use $H^{k} := W^{k,2}$ with the norm $\| \cdot \|_{H^{k}}$, and the notation $H^{1}_{0}$ and $H^{-1}$ to denote the spaces $H^{1}_{0}(\Omega)$ and its dual $H^{-1}(\Omega)$.  The duality pairing between $H^{1}_{0}$ and $H^{-1}$ is denoted by $\langle \cdot, \cdot \rangle$.

\subsection*{Useful preliminaries} We recall the Poincar\'{e} inequality for $H^{1}_{0}$:  There exists a constant $C_{p} > 0$ that depends only on $\Omega$ such that

\ \vspace*{-10pt}
\begin{align}\label{PoincareH10}
\| f \|_{L^{2}} \leq C_{p} \|\nabla f\|_{L^{2}} =: C_{p} |f|_{H^{1}_{0}} \quad \forall f \in H^{1}_{0}.
\end{align}
We also state the Sobolev embedding $H^{1} \subset L^{r}$, $r \in [1,6]$, for dimension $d = 3$:  There exists a constant $C_{s} > 0$ depending on $\Omega$ and $r$ such that
\begin{align*}
\|f\|_{L^{r}} \leq C_{s} \|f\|_{H^{1}}.
\end{align*}
We will also use the following Gronwall inequality in integral form (see \cite[Lem 3.1]{GLNeu} for a proof): Let $\alpha, \beta, u$ and $v$ be real-valued functions defined on $[0,T]$.  Assume that $\alpha$ is integrable, $\beta$ is non-negative and continuous, $u$ is continuous, $v$ is non-negative and integrable.  If $u$ and $v$ satisfy the integral inequality
\begin{align*}
u(s) + \int_{0}^{s} v(t) \, dt  \leq \alpha(s) + \int_{0}^{s} \beta(t) u(t) \, dt \quad \text{ for } s \in (0,T],
\end{align*}
then it holds that
\begin{align}
\label{Gronwall}
u(s) + \int_{0}^{s} v(t) \, dt \leq \alpha(s) + \int_{0}^{s} \beta(t) \alpha(t) \exp \left ( \int_{0}^{t} \beta(r) \, dr \right ) \, dt.
\end{align}

\section{Main results}\label{sec:mainresults}
In this section we state the main results on existence, regularity, uniqueness and continuous dependence first for regular potentials and then for singular potentials.  We also state the results for the quasi-static limit $\kappa \to 0$ for both cases.  The results are stated for dimension $d = 3$, but similar results also holds for $d = 1, 2$.

\subsection{Regular potentials}

\begin{assump}\label{assump:general}
We make the following assumptions.
\begin{enumerate}
\item[(A1)] $\lambda_{p}$, $\lambda_{a}$, $\lambda_{c}$, $\chi$ and $\eta$ are fixed non-negative constants, while $\kappa$, $\gamma$ and $\varepsilon$ are fixed positive constants.
\item[(A2)] The functions $D, h$ belong to the space $C^{0}(\mathbb{R}) \cap L^{\infty}(\mathbb{R})$, and there exist positive constants $D_{0}$, $D_{1}$ and $h_{\infty}$ such that
\begin{align*}
D_{0} \leq D(y) \leq D_{1}, \quad 0 \leq h(y) \leq h_{\infty} \quad \forall y \in \mathbb{R}.
\end{align*}
\item[(A3)] The potential $\Psi \in C^{1}(\mathbb{R})$ is a quadratic perturbation of a convex function, i.e., $\Psi(s) = \Psi_{1}(s) + \Psi_{2}(s)$ with a convex function $\Psi_{1}(s)$ and $\| \Psi_{2}'' \|_{L^{\infty}(\mathbb{R})} < \infty$.  Furthermore, $\Psi$ is non-negative and satisfies
\begin{align}\label{assump:Psi}
 |\Psi'(y)|^{s} \leq k_{1} (1 + \Psi(y)) \quad \forall y \in \mathbb{R},
\end{align}
for some constant $k_{1} > 0$, and exponent $s \in (1,2]$.
\item[(A4)] The initial and boundary data satisfy
\begin{align}
\varphi_{0} & \in H^{1}, \quad \sigma_{0} \in L^{2}, \quad \Psi(\varphi_{0}) \in L^{1}, \label{assump:initialdata} \\
\mu_{\infty}, \sigma_{\infty} & \in L^{2}(0,T;H^{1}) \cap H^{1}(0,T;L^{2}) \subset C^{0}([0,T];L^{2}). \label{assump:boundarydata}
\end{align}
\end{enumerate}
\end{assump}
We point out that by \eqref{assump:Psi} the potential $\Psi$ can have arbitrary polynomial growth (cf. \cite[\S 3 (H5)]{ColliFG}, \cite[\S 3, (H5)]{FGRNS}).

\begin{definition}\label{defn:truncated:weak}
We call a triplet of functions $(\varphi, \mu, \sigma)$ a weak solution to \eqref{Intro:CHN} if
\begin{align*}
\varphi & \in \left (-1 + L^{\infty}(0,T;H^{1}_{0}) \right ) \cap H^{1}(0,T;H^{-1}), \\
\mu & \in \mu_{\infty} + L^{2}(0,T;H^{1}_{0}), \\
\sigma & \in \left ( \sigma_{\infty} + L^{2}(0,T;H^{1}_{0}) \right ) \cap H^{1}(0,T;H^{-1}),
\end{align*}
such that $\varphi(0) = \varphi_{0}$ and $\sigma(0) = \sigma_{0}$, and satisfies for all $\zeta, \lambda, \xi \in H^{1}_{0}$ and for a.e. $t \in (0,T)$,
\begin{subequations}\label{truncated:weak}
\begin{align}
0 & = \langle \partial_{t}\varphi, \zeta \rangle + \int_{\Omega} \nabla \mu \cdot \nabla \zeta - (\lambda_{p} \sigma - \lambda_{a}) h(\varphi) \zeta \, dx, \label{truncated:weak:varphi} \\
0 & = \int_{\Omega} \mu \lambda - \frac{\gamma}{\varepsilon} \Psi'(\varphi) \lambda - \gamma \varepsilon \nabla \varphi \cdot \nabla \lambda + \chi \sigma \lambda \,dx,  \label{truncated:weak:mu} \\
0 & = \langle \kappa \partial_{t} \sigma, \xi \rangle + \int_{\Omega} D(\varphi) (\nabla \sigma - \eta \nabla \varphi) \cdot \nabla \xi + \lambda_{c} \sigma h(\varphi) \xi \, dx. \label{truncated:weak:sigma}
\end{align}
\end{subequations}
\end{definition}

\begin{theorem}[Existence of weak solutions and energy inequality]\label{thm:truncated:exist}
Let $\Omega \subset \mathbb{R}^{3}$ denote a bounded domain with $C^{2}$-boundary $\Gamma$ and suppose that Assumption \ref{assump:general} holds.  Then there exists a weak solution $(\varphi, \mu, \sigma)$ to \eqref{Intro:CHN} in the sense of Definition \ref{defn:truncated:weak} which satisfies
\begin{equation}\label{energyineq}
\begin{aligned}
& \sup_{s \in (0,T]} \left ( \| \Psi(\varphi(s)) \|_{L^{1}} + \| \varphi(s) \|_{H^{1}}^{2} + \kappa \| \sigma(s)\|_{L^{2}}^{2} \right ) \\
& \qquad + \| \mu\|_{L^{2}(0,T;H^{1})}^{2} + \|\sigma\|_{L^{2}(0,T;H^{1})}^{2} \\
& \quad \leq  C \left ( 1 + \kappa \|\sigma_{0} - \sigma_{\infty}(0)\|_{L^{2}}^{2} + \kappa^{2} \| \partial_{t} \sigma_{\infty} \|_{L^{2}(0,T;L^{2})}^{2} + \kappa \| \sigma_{\infty}\|_{L^{\infty}(0,T;L^{2})}^{2} \right ),
\end{aligned}
\end{equation}
for some positive constant $C$ independent of $\kappa$, $\varphi$, $\mu$ and $\sigma$.
\end{theorem}

\begin{theorem}[Regularity]\label{thm:truncated:regularity}
Suppose $(\varphi, \mu, \sigma)$ is a weak solution triplet to \eqref{Intro:CHN} in the sense of Definition \ref{defn:truncated:weak}.
\begin{itemize}
\item \label{Reg1} If, instead of $\mathrm{(A3)}$, $\Psi \in C^{1}(\mathbb{R})$ is non-negative and satisfies
\begin{align}\label{Psi'5thorder}
|\Psi'(y)| \leq k_{2} \left ( 1 + |y|^{p} \right ) \quad \forall y \in \mathbb{R},
\end{align}
for some constant $k_{2} > 0$ and exponent $p \in (1,5)$.  Then $\varphi \in L^{2}(0,T;W^{2,6})$.
\item \label{Reg2} If $\Gamma$ is $C^{3}$, and instead of $\mathrm{(A3)}$, $\Psi \in C^{1,1}(\mathbb{R})$ is non-negative and satisfies
\begin{align}\label{assump:Psi''}
|\Psi''(y)| \leq k_{3} (1 + |y|^{p-1}) \text{ for a.e. } y \in \mathbb{R},
\end{align}
for some constant $k_{3} > 0$ and exponent $p \in (1,5)$. Then $\varphi \in L^{2}(0,T;H^{3})$.
\end{itemize}
\end{theorem}
Note that in the above assumptions, we do not require $\Psi$ to be a quadratic perturbation of a convex function, as in $\mathrm{(A3)}$.  In fact, $\Psi$ can be a general function with the assumed regularity and polynomial growth.  We make the following additional assumptions to prove continuous dependence on the initial and boundary data.
\begin{assump}\label{assump:Ctsdep}
In addition to Assumption \ref{assump:general}, we assume that
\begin{enumerate}
\item[(C1)] $D(\cdot) = D_{0} > 0$ is a constant.
\item[(C2)] $h(\cdot)$ is Lipschitz continuous with Lipschitz constant $\mathrm{L}_{h}$.
\item[(C3)] $\Psi'$ satisfies
\begin{align}
\label{Psi'':uniqueness}
|\Psi'(s_{1}) - \Psi'(s_{2})| \leq k_{3} (1 + |s_{1}|^{4} + |s_{2}|^{4})|s_{1} - s_{2}|,
\end{align}
for some positive constant $k_{3}$ and for all $s_{1}, s_{2} \in \mathbb{R}$.
\end{enumerate}
\end{assump}

\begin{theorem}[Continuous dependence]
\label{thm:truncated:ctsdep}
Under Assumption \ref{assump:Ctsdep}, for any two weak solution triplets $\{\varphi_{i}, \mu_{i}, \sigma_{i}\}_{i = 1,2}$ to \eqref{Intro:CHN} satisfying
\begin{align*}
\varphi_{i} & \in (-1 + L^{\infty}(0,T;H^{1}_{0})) \cap L^{2}(0,T;H^{3}) \cap H^{1}(0,T;H^{-1}), \\
\mu_{i} & \in \mu_{\infty,i} + L^{2}(0,T;H^{1}_{0}), \\
\sigma_{i} & \in (\sigma_{\infty,i} + L^{2}(0,T;H^{1}_{0})) \cap L^{\infty}(0,T;L^{2}) \cap H^{1}(0,T;H^{-1}),
\end{align*}
with $\mu_{\infty,i}, \sigma_{\infty,i}$ satisfying \eqref{assump:boundarydata}, $\varphi_{i}(0) = \varphi_{i,0} \in H^{1}$ and $\sigma_{i}(0) = \sigma_{i,0} \in L^{2}$ for $i = 1,2$, there exists a positive constant $C$, depending on $T$, $\gamma$, $\varepsilon$, $D_{0}$, $\eta$, $\chi$, $\lambda_{p}$, $\lambda_{a}$, $\lambda_{c}$, $k_{3}$, $\|\varphi_{i}\|_{L^{2}(0,T;H^{3})}$, $\|\varphi_{i}\|_{L^{\infty}(0,T;H^{1})}$, $\|\sigma_{i}\|_{L^{\infty}(0,T;L^{2})}$, $\mathrm{L}_{h}$, $\Omega$, and $h_{\infty}$ such that
\begin{align*}
& \sup_{t \in (0,T]} \left ( \|\varphi_{1}(t) - \varphi_{2}(t)\|_{L^{2}}^{2} + \kappa \|\sigma_{1}(t) - \sigma_{2}(t)\|_{L^{2}}^{2} \right ) \\
& \qquad  + \|\mu_{1} - \mu_{2}\|_{L^{2}(0,T;L^{2})}^{2} + \|\sigma_{1} - \sigma_{2}\|_{L^{2}(0,T;H^{1})}^{2} + \|\varphi_{1} - \varphi_{2}\|_{L^{2}(0,T;H^{1})}^{2} \\
& \quad \leq  C \left ( \|\varphi_{1,0} - \varphi_{2,0}\|_{L^{2}}^{2} + \|\mu_{\infty,1} - \mu_{\infty,2}\|_{L^{2}(0,T;H^{1})}^{2} + \|\sigma_{\infty,1} - \sigma_{\infty,2}\|_{L^{2}(0,T;H^{1})}^{2} \right ) \\
& \qquad + C \kappa \left ( \|\sigma_{\infty,1} - \sigma_{\infty,2}\|_{L^{\infty}(0,T;L^{2})}^{2} + \|\sigma_{1,0} - \sigma_{2,0}\|_{L^{2}}^{2} \right ) \\
& \qquad + C \kappa^{2} \|\partial_{t}(\sigma_{\infty,1} - \sigma_{\infty,2})\|_{L^{2}(0,T;L^{2})}^{2}.
\end{align*}
\end{theorem}
Here, we point out that, the constant $C$ does not depend on $\kappa$, and Theorem \ref{thm:truncated:ctsdep} provides continuous dependence for the chemical potential $\mu$ in $L^{2}(0,T;L^{2})$ and for the order parameter $\varphi$ in $L^{\infty}(0,T;L^{2})$.  This is in contrast to the classical continuous dependence for $\varphi$ in $L^{\infty}(0,T;H^{-1})$ one expects for the Cahn--Hilliard equation, compare \cite[Thm. 2]{FGR} and Theorem \ref{thm:singular:ctsdep} below.

Before we give the result concerning the quasi-static limit $\kappa \to 0$, we introduce the definition of weak solutions to \eqref{Intro:quasistatic}.
\begin{definition}\label{defn:quasi:weak}
We call a triplet of functions $(\varphi_{*}, \mu_{*}, \sigma_{*})$ a weak solution to \eqref{Intro:quasistatic} if
\begin{align*}
\varphi_{*} & \in \left ( -1 + L^{\infty}(0,T;H^{1}_{0}) \right ) \cap H^{1}(0,T;H^{-1}), \\
\mu_{*} & \in \mu_{\infty} + L^{2}(0,T;H^{1}_{0}), \\
\sigma_{*} & \in \sigma_{\infty} + L^{2}(0,T;H^{1}_{0}),
\end{align*}
such that $\varphi_{*}(0) = \varphi_{0}$, and satisfies for all $\zeta, \lambda, \xi \in H^{1}_{0}$ and for a.e. $t \in (0,T)$,
\begin{subequations}
\label{quasistatic:weak}
\begin{align}
0 & = \langle \partial_{t}\varphi_{*}, \zeta \rangle + \int_{\Omega} \nabla \mu_{*} \cdot \nabla \zeta - (\lambda_{p} \sigma_{*} - \lambda_{a}) h(\varphi_{*}) \zeta \, dx, \label{quasistatic:weak:varphi} \\
0 & = \int_{\Omega} \mu_{*} \lambda - \frac{\gamma}{\varepsilon} \Psi'(\varphi_{*}) \lambda - \gamma \varepsilon \nabla \varphi_{*} \cdot \nabla \lambda + \chi \sigma_{*} \lambda \, dx,  \label{quasistatic:weak:mu} \\
0 & = \int_{\Omega} D(\varphi_{*}) (\nabla \sigma_{*} - \eta \nabla \varphi_{*}) \cdot \nabla \xi + \lambda_{c} \sigma_{*} h(\varphi_{*}) \xi \, dx. \label{quasistatic:weak:sigma}
\end{align}
\end{subequations}
\end{definition}

\begin{theorem}[Quasi-static limit]\label{thm:truncated:quasi}
For each $\kappa \in (0,1]$, let $(\varphi^{\kappa}, \mu^{\kappa}, \sigma^{\kappa})$ denote a weak solution triplet to \eqref{Intro:CHN} with initial conditions $\varphi_{0}$ and $\sigma_{0}$ satisfying \eqref{assump:initialdata}.  Then, as $\kappa \to 0$, we have
\begin{alignat*}{3}
\varphi^{\kappa} & \rightarrow \varphi_{*} && \quad \text{ weakly-}* && \quad \text{ in } L^{\infty}(0,T;H^{1}) \cap H^{1}(0,T;H^{-1}), \\
\sigma^{\kappa} & \rightarrow \sigma_{*} && \quad \text{ weakly } && \quad \text{ in } L^{2}(0,T;H^{1}), \\
\partial_{t}(\kappa \sigma^{\kappa}) & \rightarrow 0 && \quad \text{ weakly } && \quad \text{ in } L^{2}(0,T;H^{-1}), \\
\mu^{\kappa} & \rightarrow \mu_{*} && \quad \text{ weakly } && \quad \text{ in } L^{2}(0,T;H^{1}),
\end{alignat*}
such that $(\varphi_{*}, \mu_{*}, \sigma_{*})$ is a weak solution triplet to \eqref{Intro:quasistatic} in the sense of Definition \ref{defn:quasi:weak}.  Furthermore the assertions of Theorem \ref{thm:truncated:regularity} also hold for $\varphi_{*}$, and if in addition $\sigma_{\infty} \in L^{\infty}(0,T;H^{1})$ holds, then we have $\sigma_{*} \in L^{\infty}(0,T;H^{1})$.
\end{theorem}

\begin{theorem}[Continuous dependence]\label{thm:truncated:quasi:ctsdep}
Under Assumption \ref{assump:Ctsdep}, for any two weak solution triplets $\{ \varphi_{i}, \mu_{i}, \sigma_{i}\}_{i = 1,2}$ to \eqref{Intro:quasistatic} satisfying
%
\begin{align*}
\varphi_{i} & \in (-1 + L^{\infty}(0,T;H^{1}_{0})) \cap L^{2}(0,T;H^{3}) \cap  H^{1}(0,T;H^{-1}), 
\end{align*}\begin{align*}\mu_{i} & \in \mu_{\infty,i} + L^{2}(0,T;H^{1}_{0}), \\
\sigma_{i} & \in (\sigma_{\infty,i} + L^{2}(0,T;H^{1}_{0})) \cap L^{\infty}(0,T;L^{2}),
\end{align*}
with $\mu_{\infty,i}$, $\sigma_{\infty,i}$ satisfying \eqref{assump:boundarydata} and $\varphi_{i}(0) = \varphi_{i,0} \in H^{1}$ for $i = 1,2$, there exists a positive constant $C$, depending on $T$, $\gamma$, $\varepsilon$, $\chi$, $D_{0}$, $\eta$, $\lambda_{p}$, $\lambda_{a}$, $\lambda_{c}$, $k_{3}$, $\|\varphi_{i}\|_{L^{2}(0,T;H^{3})}$, $\|\varphi_{i}\|_{L^{\infty}(0,T;H^{1})}$, $\|\sigma_{i}\|_{L^{\infty}(0,T;L^{2})}$, $\mathrm{L}_{h}$, $\Omega$, and $h_{\infty}$ such that
\begin{align*}
& \sup_{t \in (0,T]} \|\varphi_{1}(t) - \varphi_{2}(t)\|_{L^{2}}^{2} + \|\mu_{1} - \mu_{2}\|_{L^{2}(0,T;L^{2})}^{2} \\
& \qquad + \|\sigma_{1} - \sigma_{2}\|_{L^{2}(0,T;H^{1})}^{2} + \|\varphi_{1} - \varphi_{2}\|_{L^{2}(0,T;H^{1})}^{2} \\
& \quad \leq C \left ( \|\varphi_{1,0} - \varphi_{2,0}\|_{L^{2}}^{2} + \|\mu_{\infty,1} - \mu_{\infty,2}\|_{L^{2}(0,T;H^{1})}^{2} + \|\sigma_{\infty,1} - \sigma_{\infty,2}\|_{L^{2}(0,T;H^{1})}^{2} \right ).
\end{align*}
\end{theorem}
By Theorem \ref{thm:truncated:quasi}, the condition $\sigma_{i} \in L^{\infty}(0,T;L^{2})$ is fulfilled if, for instance, $\sigma_{\infty,i} \in L^{\infty}(0,T;H^{1})$.

\subsection{Singular potentials}

\begin{definition}\label{defn:singular:weak}
We call a quadruple of functions $(\varphi, \mu, \sigma, \psi)$ a weak solution to \eqref{singular:CHN} if
\begin{align*}
\varphi & \in \left ( -1 + L^{\infty}(0,T;H^{1}_{0}) \right ) \cap H^{1}(0,T;H^{-1}), \\
\mu & \in \mu_{\infty} + L^{2}(0,T;H^{1}_{0}),  \\
\sigma & \in \left ( \sigma_{\infty} + L^{2}(0,T;H^{1}_{0})  \right ) \cap H^{1}(0,T;H^{-1}), \\
\psi & \in L^{2}(0,T;L^{2}) \text{ with } \varphi \in D(\beta), \psi \in \beta(\varphi) \text{ a.e. in } \Omega \times (0,T),
\end{align*}
such that $\varphi(0) = \varphi_{0}$ and $\sigma(0) = \sigma_{0}$, and satisfies for all $\zeta, \lambda, \xi \in H^{1}_{0}$ and for a.e. $t \in (0,T)$,
\begin{subequations}\label{singular:weak}
\begin{align}
0 & = \langle \partial_{t}\varphi, \zeta \rangle + \int_{\Omega} \nabla \mu \cdot \nabla \zeta - (\lambda_{p} \sigma - \lambda_{a}) h(\varphi) \zeta \, dx, \label{singular:weak:varphi} \\
0 & = \int_{\Omega} \mu \lambda - \frac{\gamma}{\varepsilon} (\psi + \Lambda'(\varphi)) \lambda - \gamma \varepsilon \nabla \varphi \cdot \nabla \lambda + \chi \sigma \lambda \, dx,  \label{singular:weak:mu} \\
0 & = \langle \kappa \partial_{t} \sigma, \xi \rangle + \int_{\Omega} D(\varphi) (\nabla \sigma - \eta \nabla \varphi) \cdot \nabla \xi + \lambda_{c} \sigma h(\varphi) \xi \, dx. \label{singular:weak:sigma}
\end{align}
\end{subequations}
\end{definition}

\begin{assump}\label{assump:singularpotential}
In addition to $(\mathrm{A1})$, $(\mathrm{A2})$, \eqref{assump:boundarydata} and $\sigma_{0} \in L^{2}$, we assume that $\Psi$ is of the form \eqref{defn:singularpot} with
\begin{enumerate}
\item[(S1)] $\hat{\beta}: \mathbb{R} \to [0,\infty]$ is a convex, proper, lower semicontinuous function  with $\hat{\beta}(0) = 0$ and $-1 \in D(\beta)$.
\item[(S2)] $\Lambda \in C^{1,1}(\mathbb{R})$ is non-negative with $\|\Lambda''\|_{L^{\infty}(\mathbb{R})} < \infty$.
\item[(S3)] $\varphi_{0} \in H^{1}$ satisfies $\hat{\beta}(\varphi_{0}) \in L^{1}$.
\end{enumerate}
\end{assump}
Note that from the definition \eqref{subdifferential}, the condition $\hat{\beta}(0) = 0$ and non-negativity of $\hat{\beta}$ imply that $0 \in \beta(0)$.  Meanwhile, the condition $-1 \in D(\beta)$ is motivated from the boundary condition for $\varphi$ and in particular this implies that $\beta(-1) \neq \emptyset$ and $|\beta^{0}(-1)| < \infty$.  This is required to obtain estimates on the selection $\psi$ in $L^{2}(0,T;L^{2})$, and while this condition is satisfied for the double-obstacle potential \eqref{DoubleObstaclebeta}, unfortunately it does not hold for the classical logarithmic potential
\begin{align*}
\hat{\beta}_{\mathrm{log}}(y) = (1-y) \log (1-y) + (1+y) \log(1+y), \quad \Lambda_{\mathrm{log}}(y) = \theta (1-y^{2}) \text{ for } \theta > 0,
\end{align*}
whose subdifferential $\beta_{\mathrm{log}}$ has an effective domain $D(\beta_{\mathrm{log}}) = (-1,1)$.  Thus, our current setting does not extend to the logarithmic potential, but this can be remedied if we impose that the more general boundary condition $\varphi_{\infty}$ for $\varphi$ lies in $D(\beta)$.

\begin{theorem}[Existence of regular weak solutions and energy inequality]\label{thm:singular:exist}
Let $\Gamma$ be a $C^{3}$-boundary.  Suppose that Assumption \ref{assump:singularpotential} is satisfied.  Then, there exists a weak solution $(\varphi, \mu, \sigma, \psi)$ to \eqref{singular:CHN} in the sense of Definition \ref{defn:singular:weak} with the additional regularities
\begin{align*}
\varphi \in L^{2}(0,T;H^{2}), \quad
\sigma \in L^{\infty}(0,T;L^{2}).
\end{align*}
Furthermore, it holds that
\begin{equation}\label{energyineq:singular}
\begin{aligned}
& \sup_{s \in (0,T]} \left (\int_{\Omega}  \hat{\beta}(\varphi(s)) + \Lambda(\varphi(s)) \, dx + \|\varphi(s)\|_{H^{1}}^{2} + \kappa \|\sigma(s)\|_{L^{2}}^{2} \right ) \\
& \qquad + \|\mu\|_{L^{2}(0,T;H^{1})}^{2} + \|\sigma\|_{L^{2}(0,T;H^{1})}^{2} \\
&\quad  \leq C \left ( 1 + \kappa \|\sigma_{0} - \sigma_{\infty}(0)\|_{L^{2}}^{2} + \kappa^{2} \|\partial_{t} \sigma_{\infty}\|_{L^{2}(0,T;L^{2})}^{2} + \kappa \|\sigma_{\infty}\|_{L^{\infty}(0,T;L^{2})}^{2} \right ),
\end{aligned}
\end{equation}
for some positive constant $C$ independent of $\kappa$, $\varphi$, $\mu$, $\sigma$ and $\psi$.
\end{theorem}

We point out that we are only able to show $H^{2}$-regularity for $\varphi$, in contrast with the $H^{3}$-regularity for the regular potentials.  This is due to the fact that the $H^{3}$-regularity estimate is not independent of the Yosida parameter.  Before we state the result on continuous dependence, we introduce the inverse Dirichlet-Laplacian operator $N: H^{-1} \to H^{1}_{0}$, $f \mapsto N(f)$, via
\begin{align}\label{InverseLaplacian}
\int_{\Omega} \nabla N(f) \cdot \nabla u \, dx = \langle f, u \rangle \quad \forall u \in H^{1}_{0}.
\end{align}
That is, $N(f)$ is the unique solution to the Dirichlet-Laplacian problem with right-hand side $f$.  Note that
\begin{align*}
\langle f, N(g) \rangle = \int_{\Omega} \nabla N(f) \cdot \nabla  N(g) \, dx = \langle g, N(f) \rangle \quad \forall f,g \in H^{-1},
\end{align*}
and we can define a norm on $H^{-1}$ as
\begin{align}\label{starnorm}
\|f\|_{*} := \|\nabla N(f)\|_{L^{2}} = |N(f)|_{H^{1}_{0}} = \sqrt{\langle f, N(f) \rangle},
\end{align}
which satisfies for $g \in H^{1}_{0}$, \begin{align}
\label{L2normstarnormgradnorm} \|g\|_{L^{2}}^{2}  = \int_{\Omega} \nabla N(g) \cdot \nabla g \, dx \leq \|\nabla N(g)\|_{L^{2}} \|\nabla g\|_{L^{2}} & = \|g\|_{*} \|\nabla g\|_{L^{2}}.
\end{align}
Moreover, we obtain from \eqref{InverseLaplacian} the relation
\begin{align}\label{TimederivativeStar}
\frac{1}{2} \frac{d}{dt} \|f(t)\|_{*}^{2} = \frac{1}{2} \frac{d}{dt} \|\nabla N(f(t))\|_{L^{2}}^{2} = \langle \partial_{t} f(t), N(f(t)) \rangle
\end{align}
for $f \in H^{1}(0,T;H^{-1})$ and a.e. $t \in (0,T)$.  This follows from the fact that $\partial_{t}f \in L^{2}(0,T;H^{-1})$ and thus $N(\partial_{t}f)$ is well-defined and satisfies
\begin{align*}
\int_{\Omega} \nabla N(\partial_{t}f) \cdot \nabla u \, dx = \langle \partial_{t}f, u \rangle \quad \forall u \in H^{1}_{0}.
\end{align*}
Multiplying the above equality with $\zeta \in C^{\infty}_{c}(0,T)$ and integrating in time, we obtain the following chain of equalities (here $(\cdot, \cdot)$ denotes the $L^{2}$-scalar product)
\begin{align*}
\int_{0}^{T} \zeta (\nabla N(\partial_{t}f),\nabla u) \, dt & = \left \langle \int_{0}^{T} \partial_{t}f \zeta \, dt, u \right \rangle = \left \langle -\int_{0}^{T} f \zeta' \, dt, u \right \rangle \\
&  = -\int_{0}^{T} \zeta' (\nabla N(f), \nabla u )\,dt.
\end{align*}
As $u \in H^{1}_{0}$ and $\zeta \in C^{\infty}_{c}(0,T)$ are arbitrary, we find that $\partial_{t}(\nabla N(f)) = \nabla N(\partial_{t}f)$ in the weak sense.  Thus
\begin{align*}
\frac{1}{2} \frac{d}{dt} \|\nabla N(f)\|_{L^{2}}^{2} = (\partial_{t}(\nabla N(f)),\nabla N(f)) = \langle \partial_{t}f, N(f) \rangle.
\end{align*}

\begin{theorem}[Partial continuous dependence and full uniqueness]\label{thm:singular:ctsdep}
Under $(\mathrm{C1})$ and $(\mathrm{C2})$ of Assumption \ref{assump:Ctsdep}, for any two weak solution quadruples $\{\varphi_{i}, \mu_{i}, \sigma_{i}, \psi_{i}\}_{i = 1,2}$ to \eqref{singular:CHN} obtained from Theorem \ref{thm:singular:exist} with boundary conditions $\mu_{\infty,1} = \mu_{\infty,2}$,\break $\{\sigma_{\infty,i}\}_{i=1,2}$ satisfying \eqref{assump:boundarydata} and initial conditions $\varphi_{i}(0) = \varphi_{i,0} \in H^{1}$ and $\sigma_{i}(0) = \sigma_{i,0} \in L^{2}$ for $i = 1,2$, there exists a positive constant $C$, depending on $T$, $\gamma$, $\varepsilon$, $\chi$, $D_{0}$, $\eta$, $\lambda_{p}$, $\lambda_{a}$, $\lambda_{c}$, $\|\sigma_{i}\|_{L^{\infty}(0,T;L^{2})}$, $\mathrm{L}_{h}$, $\Omega$, $\|\Lambda''\|_{L^{\infty}(\mathbb{R})}$, and $h_{\infty}$ such that
\begin{align*}
& \sup_{t \in (0,T]} \left ( \|\varphi_{1}(t) - \varphi_{2}(t)\|_{*}^{2} + \kappa \|\sigma_{1}(t) - \sigma_{2}(t)\|_{L^{2}}^{2} \right ) \\
& \qquad + \|\sigma_{1} - \sigma_{2}\|_{L^{2}(0,T;H^{1})}^{2} + \|\varphi_{1} - \varphi_{2}\|_{L^{2}(0,T;H^{1})}^{2} \\
& \quad \leq C \left ( \|\varphi_{1,0} - \varphi_{2,0}\|_{*}^{2} + \|\sigma_{\infty,1} - \sigma_{\infty,2}\|_{L^{2}(0,T;H^{1})}^{2} + \kappa^{2}  \|\sigma_{\infty,1} - \sigma_{\infty,2}\|_{H^{1}(0,T;L^{2})}^{2}\right ) \\
& \qquad + C \kappa \left ( \|\sigma_{\infty,1} - \sigma_{\infty,2}\|_{L^{\infty}(0,T;L^{2})}^{2} + \|\sigma_{1,0} - \sigma_{2,0}\|_{L^{2}}^{2}  \right ).
\end{align*}
In particular, if $\sigma_{\infty,1} = \sigma_{\infty,2}$, $\varphi_{1,0} = \varphi_{2,0}$ and $\sigma_{1,0} = \sigma_{2,0}$ hold, then we have $\varphi_{1} = \varphi_{2}$, $\mu_{1} = \mu_{2}$, $\sigma_{1} = \sigma_{2}$ and $\psi_{1} = \psi_{2}$.
\end{theorem}

We note that continuous dependence on the initial and boundary data can only be shown for $\varphi$ and $\sigma$.  In particular, we do not have control over the differences $\mu_{1} - \mu_{2}$ and $\psi_{1} - \psi_{2}$.  In the proof we have to apply the inverse Dirichlet-Laplacian operator \eqref{InverseLaplacian} and thus it is necessary that $\mu_{\infty,1} = \mu_{\infty,2}$.  Partial continuous dependence are often observed in the case $\beta$ is multivalued (see for instance  \cite[Remark 2.3]{Gilardi} and \cite[Remark 3.1]{GMS}).

In the same spirit as Theorem \ref{thm:truncated:quasi}, we will pass to the limit $\kappa \to 0$ in \eqref{singular:CHN} to deduce the existence of weak solutions for the quasi-static model of \eqref{singular:CHN} with $\kappa = 0$.

\begin{theorem}[Quasi-static limit with singular potentials]\label{thm:singular:quasi}
\
\begin{itemize}
\item $(\textrm{Existence})$
For each $\kappa \in (0,1]$, let $(\varphi^{\kappa}, \mu^{\kappa}, \sigma^{\kappa}, \psi^{\kappa})$ denote a weak solution to \eqref{singular:CHN} obtained from Theorem \ref{thm:singular:exist} with initial conditions $\sigma_{0} \in L^{2}$, $\varphi_{0} \in H^{1}$ such that $\hat{\beta}(\varphi_{0}) \in L^{1}$.  Then, as $\kappa \to 0$, we have
\begin{alignat*}{3}
\varphi^{\kappa} & \rightarrow \overline{\varphi} && \quad \text{ weakly-}* && \quad \text{ in } L^{\infty}(0,T;H^{1}) \cap H^{1}(0,T;H^{-1}) \cap L^{2}(0,T;H^{2}), \\
\sigma^{\kappa} & \rightarrow \overline{\sigma} && \quad \text{ weakly } && \quad \text{ in } L^{2}(0,T;H^{1}), \\
\partial_{t}(\kappa \sigma^{\kappa}) & \rightarrow 0 && \quad \text{ weakly } && \quad \text{ in } L^{2}(0,T;H^{-1}), \\
\mu^{\kappa} & \rightarrow \overline{\mu} && \quad \text{ weakly } && \quad \text{ in } L^{2}(0,T;H^{1}), \\
\psi^{\kappa} & \rightarrow \overline{\psi} && \quad \text{ weakly } && \quad \text{ in } L^{2}(0,T;L^{2}),
\end{alignat*}
such that $(\overline{\varphi}, \overline{\mu}, \overline{\sigma}, \overline{\psi})$ satisfies \eqref{singular:weak} with $\kappa = 0$.
\vspace{2mm}
\item $(\textrm{Regularity})$ If $\sigma_{\infty} \in L^{\infty}(0,T;H^{1})$ then $\overline{\sigma} \in L^{\infty}(0,T;H^{1})$.
\vspace{2mm}
\item $(\textrm{Partial continuous dependence})$  Under $(\mathrm{C1})$ and $(\mathrm{C2})$ of Assumption \ref{assump:Ctsdep}, for any two weak solution quadruples $\{\overline{\varphi}_{i}, \overline{\mu}_{i}, \overline{\sigma}_{i}, \overline{\psi}_{i}\}_{i=1,2}$ to \eqref{singular:CHN} with $\kappa = 0$ satisfying additionally $\overline{\sigma}_{i} \in L^{\infty}(0,T;L^{2})$, the corresponding boundary conditions $\mu_{\infty,1}= \mu_{\infty,2}$, $\{\sigma_{\infty,i}\}_{i=1,2}$ and initial conditions $\varphi_{i,0} \in H^{1}$ for $i = 1,2$, there exists a positive constant $C$, depending on $T$, $\gamma$, $\varepsilon$, $\chi$, $D_{0}$, $\eta$, $\lambda_{p}$, $\lambda_{a}$, $\lambda_{c}$, $\|\sigma_{i}\|_{L^{\infty}(0,T;L^{2})}$, $\mathrm{L}_{h}$, $\Omega$, $\|\Lambda''\|_{L^{\infty}(\mathbb{R})}$, and $h_{\infty}$ such that
\begin{align*}
& \sup_{t \in (0,T]}  \|\varphi_{1}(t) - \varphi_{2}(t)\|_{*}^{2} + \|\sigma_{1} - \sigma_{2}\|_{L^{2}(0,T;H^{1})}^{2} + \|\varphi_{1} - \varphi_{2}\|_{L^{2}(0,T;H^{1})}^{2} \\
& \quad \leq C \left ( \|\varphi_{1,0} - \varphi_{2,0}\|_{*}^{2} + \|\sigma_{\infty,1} - \sigma_{\infty,2}\|_{L^{2}(0,T;H^{1})}^{2} \right ).
\end{align*}
\item $(\textrm{Uniqueness})$ In particular, if $\sigma_{\infty,1} = \sigma_{\infty,2}$, and $\varphi_{1,0} = \varphi_{2,0}$ hold, then $\overline{\varphi_{1}} = \overline{\varphi_{2}}$, $\overline{\mu_{1}} = \overline{\mu_{2}}$, $\overline{\sigma_{1}} = \overline{\sigma_{2}}$ and $\overline{\psi_{1}} = \overline{\psi_{2}}$.
\end{itemize}
\end{theorem}
Since the above result is proved by combining the proofs of Theorems \ref{thm:truncated:quasi}, \ref{thm:truncated:quasi:ctsdep} and \ref{thm:singular:ctsdep}, we will omit the proof of Theorem \ref{thm:singular:quasi}.

\section{Regular potentials}\label{sec:regularPot}
\subsection{Galerkin approximation}
We prove Theorem \ref{thm:truncated:exist} via a Galerkin procedure, assuming first that $\varphi_{0} + 1 \in H^{2} \cap H^{1}_{0}$.  The general case $\varphi_{0} \in H^{1}$ with $\Psi(\varphi_{0}) \in L^{1}$ can be handled by means of a density argument, see for example \cite[Remark 2]{ColliFG}, as by $\mathrm{(A3)}$ $\Psi$ is a quadratic perturbation of a convex function.  Let us consider $\{w_{j}\}_{j \in \mathbb{N}}$ as the set of eigenfunctions of the Dirichlet-Laplacian which are chosen such that they form an orthonormal basis of $L^{2}$ and an orthogonal basis of $H^{1}_{0}$.  Furthermore, as $\Gamma$ is $C^{2}$, elliptic regularity yields that $w_{j} \in H^{2} \cap H^{1}_{0}$ for all $j \in \mathbb{N}$.  Let $W_{k} = \mathrm{span}\{w_{1}, \dots, w_{k}\} \subset H^{1}_{0}$ denote the finite dimensional subspace spanned by the first $k$ eigenfunctions.  We introduce the Galerkin ansatz
\begin{align*}
\varphi_{k} := -1 + \sum_{i=1}^{k} \alpha_{i}^{k} w_{i}, \quad \mu_{k} := \mu_{\infty} + \sum_{i=1}^{k} \beta_{i}^{k} w_{i}, \quad \sigma_{k} := \sigma_{\infty} +  \sum_{i=1}^{k} \tau_{i}^{k} w_{i},
\end{align*}
and we require that $\varphi_{k}, \mu_{k}$ and $\sigma_{k}$ satisfy
\begin{subequations}\label{Truncated:discrete}
\begin{align}
0 & = \int_{\Omega} \partial_{t}\varphi_{k} w_{j} + \nabla \mu_{k} \cdot \nabla w_{j} - (\lambda_{p} \sigma_{k} - \lambda_{a}) h(\varphi_{k}) w_{j} \, dx, \label{discrete:varphi} \\
0 & = \int_{\Omega} \mu_{k} w_{j} - \frac{\gamma}{\varepsilon} \Psi'(\varphi_{k}) w_{j} - \gamma \varepsilon \nabla \varphi_{k} \cdot \nabla w_{j} + \chi \sigma_{k} w_{j} \, dx, \label{discrete:mu} \\
0 & = \int_{\Omega} \kappa \partial_{t} \sigma_{k} w_{j} + D(\varphi_{k}) (\nabla \sigma_{k} - \eta \nabla \varphi_{k}) \cdot \nabla w_{j} + \lambda_{c} \sigma_{k} h(\varphi_{k}) w_{j} \,dx, \label{discrete:sigma}
\end{align}
\end{subequations}
for all $1 \leq j \leq k$. For convenience we introduce the following matrices
\begin{alignat*}{3}
(\boldsymbol{M})_{ij} &:= \int_{\Omega} w_{i} w_{j} \, dx = \delta_{ij}, && \quad (\boldsymbol{S})_{ij} && := \int_{\Omega} \nabla w_{i} \cdot \nabla w_{j} \, dx, \\
(\boldsymbol{M}_{h})_{ij} & := \int_{\Omega} h(\varphi_{k}) w_{i} w_{j} \, dx, && \quad (\boldsymbol{S}_{D})_{ij} && := \int_{\Omega} D(\varphi_{k}) \nabla w_{i} \cdot \nabla w_{j} \, dx,
\end{alignat*}
where $\delta_{ij}$ is the Kronecker delta.  Note that $\boldsymbol{M}$ is the identity matrix precisely due to the orthonormality of $\{w_{j}\}_{j \in \mathbb{N}}$ in $L^{2}$.  Furthermore, we use the notation
\begin{align*}
\boldsymbol{\alpha}^{k} := (\alpha_{i}^{k})_{1 \leq i \leq k}, \quad \boldsymbol{\beta}^{k} := (\beta_{i}^{k})_{1 \leq i \leq k}, \quad \boldsymbol{\tau}^{k} := (\tau_{i}^{k})_{1 \leq i \leq k},
\end{align*}
and
\begin{alignat*}{3}
\boldsymbol{\psi}^{k} & := \left ( \int_{\Omega} \Psi'(\varphi_{k}) w_{j} \, dx \right )_{1 \leq j \leq k}, && \quad \boldsymbol{h}^{k} && := \left ( \int_{\Omega} h(\varphi_{k}) w_{j} \, dx \right)_{1 \leq j \leq k}, \\
\boldsymbol{s}_{\mu} & := \left ( \int_{\Omega} \nabla \mu_{\infty} \cdot \nabla w_{j} \, dx \right )_{1 \leq j \leq k}, && \quad \boldsymbol{m}_{\mu} && := \left ( \int_{\Omega} \mu_{\infty} w_{j} \, dx \right )_{1 \leq j \leq k},\\
\boldsymbol{m}_{\sigma} & := \left ( \int_{\Omega} \sigma_{\infty} w_{j} \, dx \right )_{1 \leq j \leq k}, && \quad \boldsymbol{m}_{h,\sigma}^{k} && := \left ( \int_{\Omega} h(\varphi_{k}) \sigma_{\infty} w_{j} \, dx \right )_{1 \leq j \leq k}, \\
\boldsymbol{m}_{\partial_{t}\sigma} & := \left ( \int_{\Omega} \partial_{t} \sigma_{\infty} w_{j} \, dx \right )_{1 \leq j \leq k}, && \quad \boldsymbol{s}_{D,\sigma}^{k} && := \left ( \int_{\Omega} D(\varphi_{k}) \nabla \sigma_{\infty} \cdot \nabla w_{j} \, dx \right )_{1 \leq j \leq k},
\end{alignat*}
so that when we substitute the Galerkin ansatz into  \eqref{Truncated:discrete}, we obtain the following nonlinear initial value problem in vector form
\begin{subequations}\label{Truncated:ODE}
\begin{alignat}{2}
\frac{d}{dt} \boldsymbol{\alpha}^{k} & = - \boldsymbol{S} \boldsymbol{\beta}^{k} + \lambda_{p} \boldsymbol{M}_{h} \boldsymbol{\tau}^{k} - \lambda_{a} \boldsymbol{h}^{k} - \boldsymbol{s}_{\mu} + \lambda_{p} \boldsymbol{m}_{h, \sigma}^{k}, \label{ODE:varphi} \\
0 & = \boldsymbol{\beta}^{k} - \frac{\gamma}{\varepsilon} \boldsymbol{\psi}^{k} - \gamma \varepsilon \boldsymbol{S} \boldsymbol{\alpha}^{k} + \chi \boldsymbol{\tau}^{k} + \boldsymbol{m}_{\mu} + \chi \boldsymbol{m}_{\sigma},  \label{ODE:mu} \\
\kappa \frac{d}{dt} \boldsymbol{\tau}^{k} & = - \boldsymbol{S}_{D} (\boldsymbol{\tau}^{k} - \eta \boldsymbol{\alpha}^{k}) - \lambda_{c} \boldsymbol{M}_{h} \boldsymbol{\tau}^{k} - \boldsymbol{s}_{D,\sigma}^{k} - \lambda_{c} \boldsymbol{m}_{h,\sigma}^{k} - \kappa \boldsymbol{m}_{\partial_{t}\sigma}, \label{ODE:sigma}
\end{alignat}
\end{subequations}
with initial data
\begin{align*}
\boldsymbol{\alpha}^{k}(0) = \left ( \int_{\Omega} (\varphi_{0} + 1) w_{j} \, dx \right)_{1 \leq j \leq k}, \quad \boldsymbol{\tau}^{k}(0) = \left ( \int_{\Omega} (\sigma_{0} - \sigma_{\infty}) w_{j} \, dx \right )_{1 \leq j \leq k},
\end{align*}
and we define $\varphi_{k}(0) := - 1 + \sum_{i=1}^{k} (\boldsymbol{\alpha}^{k}(0))_{i} w_{i}$, $\sigma_{k}(0) := \sigma_{\infty} + \sum_{i=1}^{k} (\boldsymbol{\tau}^{k}(0))_{i} w_{i}$.

Note that $\boldsymbol{\psi}^{k}$, $\boldsymbol{h}^{k}$, $\boldsymbol{m}_{h,\sigma}^{k}$, $\boldsymbol{s}_{D,\sigma}^{k}$, $\boldsymbol{M}_{h}$ and $\boldsymbol{S}_{D}$ depend nonlinearly on the solution.  We can also express \eqref{Truncated:ODE} as an initial value problem just in terms of $\boldsymbol{\alpha}^{k}$ and $\boldsymbol{\tau}^{k}$ by substituting \eqref{ODE:mu} into \eqref{ODE:varphi}.  Moreover, by the continuity of $\Psi(\cdot)$, $h(\cdot)$ and $D(\cdot)$, we observe that the right-hand side of \eqref{ODE:varphi} and \eqref{ODE:sigma} depends continuously on $\boldsymbol{\alpha}^{k}$ and $\boldsymbol{\tau}^{k}$.  Thus, by the Cauchy--Peano theorem we infer the existence of local solutions $\boldsymbol{\alpha}^{k}, \boldsymbol{\beta}^{k}, \boldsymbol{\tau}^{k} \in (C^{1}([0,t_{k}]))^{k}$ defined on a time interval $[0, t_{k}]$ for each $k \in \mathbb{N}$.

Next, we derive some a priori estimates to deduce that $\boldsymbol{\alpha}^{k}$, $\boldsymbol{\beta}^{k}$ and $\boldsymbol{\tau}^{k}$ can be extended to the interval $[0,T]$ and to allow us to pass to the limit.  In the following, the various constants $C$ may vary from line to line, but they do not depend on $k$ or on $\kappa$.

\subsubsection{A priori estimates}\label{sec:aprioriest}
We define a new variable
\begin{align*}
\omega_{k} := \sigma_{k} - \sigma_{\infty},
\end{align*}
and then multiple \eqref{discrete:varphi} with $\beta_{i}^{k} + \chi \tau_{i}^{k}$, \eqref{discrete:mu} with $\frac{d}{dt}\alpha_{i}^{k}$ and \eqref{discrete:sigma} with $\mathcal{X} \tau_{i}^{k}$ for some positive constant $\mathcal{X}$ yet to be determined.  Sum from $i = 1$ to $k$, and summing the three equations, integrating in time from $0$ to $s \in (0,T]$ leads to
\begin{align*}
&  \int_{\Omega} \left [ \frac{\gamma}{\varepsilon} \Psi(\varphi_{k}) + \frac{\gamma \varepsilon}{2} |\nabla \varphi_{k}|^{2} + \mathcal{X} \frac{\kappa}{2} |\omega_{k}|^{2} \right ](s) \, dx  \\
& \quad + \|\nabla \mu_{k}\|_{L^{2}(0,s;L^{2})}^{2} + \mathcal{X} D_{0} \|\nabla \omega_{k}\|_{L^{2}(0,s;L^{2})}^{2} + \lambda_{c} \int_{0}^{s} \int_{\Omega} h(\varphi_{k}) |\omega_{k}|^{2} \, dx \, dt
\end{align*}
\begin{equation}\label{discrete:apriori:integrated}
\begin{aligned}
& \leq \int_{\Omega} \varphi_{k}(s) (\mu_{\infty} + \chi \sigma_{\infty})(s) - \varphi_{k}(0) (\mu_{\infty} + \chi \sigma_{\infty})(0) \, dx  \\
& \qquad - \int_{0}^{s} \int_{\Omega} \varphi_{k} \partial_{t}(\mu_{\infty} + \chi \sigma_{\infty}) \, dx \, dt  \\
& \quad + \int_{0}^{s} \int_{\Omega} - \chi \nabla \mu_{k} \cdot \nabla \omega_{k} + \nabla \mu_{k} \cdot \nabla \mu_{\infty} \, dx \, dt \\
& \qquad + \int_{0}^{s} \int_{\Omega} h(\varphi_{k}) (\lambda_{p} \sigma_{k} - \lambda_{a}) (\mu_{k} - \mu_{\infty} + \chi \omega_{k})   \, dx \, dt  \\
& \quad + \mathcal{X} \int_{0}^{s} \int_{\Omega}  D(\varphi_{k}) \eta \nabla \varphi_{k} \cdot \nabla \omega_{k}  - \kappa \partial_{t} \sigma_{\infty} \omega_{k} \, dx \, dt \\
& \qquad - \mathcal{X} \int_{0}^{s} \int_{\Omega} D(\varphi_{k}) \nabla \sigma_{\infty} \cdot \nabla \omega_{k} + \lambda_{c} h(\varphi_{k}) \sigma_{\infty} \omega_{k} \, dx \, dt \\
& \quad + \int_{\Omega} \left [ \frac{\gamma}{\varepsilon} \Psi(\varphi_{k}) + \frac{\gamma \varepsilon}{2} |\nabla \varphi_{k}|^{2} + \mathcal{X} \frac{\kappa}{2} |\sigma_{k} - \sigma_{\infty}|^{2} \right ] (0) \, dx \\
& =: I_{1a} + I_{1b} + I_{2a} + I_{2b} + \mathcal{X} (I_{3a} + I_{3b}) + I_{4}.
\end{aligned}
\end{equation}

Since $\varphi_{0} + 1 \in H^{2} \cap H^{1}_{0}$, via a similar argument employed in \cite[\S 3]{GLDarcy}, one can show that $\Delta (\varphi_{k}(0) + 1)$ converges strongly to $\Delta (\varphi_{0} + 1)$ in $L^{2}$, and hence $\varphi_{k}(0) \to \varphi_{0}$ strongly in $H^{2} \subset L^{\infty}$.  Then, there exists a constant $C > 0$ such that $\| \Psi(\varphi_{k}(0)) \|_{L^{\infty}} \leq C$ for all $k \in \mathbb{N}$.  Together with \eqref{assump:initialdata} and $\sigma_{\infty} \in C^{0}([0,T];L^{2})$, we see that
\begin{align}\label{discrete:apriori:RHS:line4}
|I_{4}| \leq C \left ( 1 + \kappa \|\sigma_{0} - \sigma_{\infty}(0)\|_{L^{2}}^{2} \right ).
\end{align}
Application of H\"{o}lder's inequality, Young's inequality and the Poincar\'{e} inequality \eqref{PoincareH10} applied to $f = \varphi_{k} + 1$ leads to
\begin{equation}\label{discrete:apriori:RHS:line1}
\begin{aligned}
|I_{1a} + I_{1b} | & \leq \frac{1}{2} \|\varphi_{k}\|_{L^{2}(0,s;L^{2})}^{2} + \frac{1}{2} \|\partial_{t}(\mu_{\infty} + \chi \sigma_{\infty})\|_{L^{2}(0,s;L^{2})}^{2} + \frac{\gamma \varepsilon}{8C_{p}^{2}}\|\varphi_{k}(s)\|_{L^{2}}^{2} \\
& \quad + \frac{2C_{p}^{2}}{\varepsilon \gamma} \|(\mu_{\infty} + \chi \sigma_{\infty})(s)\|_{L^{2}}^{2} + \|\varphi_{k}(0)\|_{L^{2}}^{2} + \frac{1}{4} \|(\mu_{\infty} + \chi \sigma_{\infty})(0)\|_{L^{2}}^{2} \\
& \leq \frac{1}{2} \|\varphi_{k}\|_{L^{2}(0,s;L^{2})}^{2} + \frac{\gamma \varepsilon}{8 C_{p}^{2}} \|\varphi_{k}(s)\|_{L^{2}}^{2} + C \\
& \leq C_{p}^{2} \|\nabla \varphi_{k}\|_{L^{2}(0,s;L^{2})}^{2} + \frac{\gamma \varepsilon}{4} \|\nabla \varphi_{k}(s)\|_{L^{2}}^{2} + C .
\end{aligned}
\end{equation}
By H\"{o}lder's inequality, we have
\begin{align*}
|I_{2a}| & \leq  \|\nabla \mu_{k}\|_{L^{2}(0,s;L^{2})} \left ( \chi \|\nabla \omega_{k}\|_{L^{2}(0,s;L^{2})} + \|\nabla \mu_{\infty}\|_{L^{2}(0,s;L^{2})} \right ), \\
|I_{2b}| &  \leq  h_{\infty} \lambda_{a} |\Omega|^{\frac{1}{2}} s^{\frac{1}{2}} \left (\|\mu_{k}\|_{L^{2}(0,s;L^{2})} + \|\mu_{\infty}\|_{L^{2}(0,s;L^{2})} + \chi \|\omega_{k}\|_{L^{2}(0,s;L^{2})} \right ) \\
& \quad +  h_{\infty} \lambda_{p} \|\sigma_{k}\|_{L^{2}(0,s;L^{2})} \left ( \|\mu_{k}\|_{L^{2}(0,s;L^{2})} + \|\mu_{\infty}\|_{L^{2}(0,s;L^{2})} + \chi \|\omega_{k}\|_{L^{2}(0,s;L^{2})} \right ) .
\end{align*}
Then, by Young's inequality and the Poincar\'{e} inequality \eqref{PoincareH10} applied to $f = \mu_{k} - \mu_{\infty}$ and $f = \omega_{k}$, we obtain
\begin{align*}
|I_{2a} + I_{2b}| & \leq \frac{4}{8} \|\nabla \mu_{k}\|_{L^{2}(0,s;L^{2})}^{2} + 2 \chi^{2} \|\nabla \omega_{k}\|_{L^{2}(0,s;L^{2})}^{2} + 2 \chi^{2} \|\omega_{k}\|_{L^{2}(0,s;L^{2})}^{2} \\
& \quad  + h_{\infty}^{2} \lambda_{p}^{2} (4C_{p}^{2} + 1) \|\sigma_{k}\|_{L^{2}(0,s;L^{2})}^{2} \\
& \quad + C \left (1 + \|\mu_{\infty}\|_{L^{2}(0,T;H^{1})}^{2} + \|\sigma_{\infty}\|_{L^{2}(0,T;H^{1})}^{2} \right )
\end{align*}
\begin{equation}\label{discrete:apriori:RHS:line2}
\begin{aligned}
& \leq  \left ( 2 \chi^{2} (1 + C_{p}^{2}) + 2 C_{p}^{2}  h_{\infty}^{2} \lambda_{p}^{2} (4C_{p}^{2} + 1) \right )  \|\nabla \omega_{k}\|_{L^{2}(0,s;L^{2})}^{2} \\
& \quad + \frac{1}{2} \|\nabla \mu_{k}\|_{L^{2}(0,s;L^{2})}^{2} + C,
\end{aligned}
\end{equation}
where we have used
\begin{align}\label{sigmakomegakrelation}
\|\sigma_{k}\|_{L^{2}} \leq \|\omega_{k}\|_{L^{2}} + \|\sigma_{\infty}\|_{L^{2}}.
\end{align}
Similarly,
\begin{equation}\label{discrete:apriori:RHS:line3}
\begin{aligned}
|I_{3a} + I_{3b}| & \leq D_{1} \left (\eta \|\nabla \varphi_{k}\|_{L^{2}(0,s;L^{2})} + \|\nabla \sigma_{\infty}\|_{L^{2}(0,s;L^{2})} \right ) \|\nabla \omega_{k}\|_{L^{2}(0,s;L^{2})} \\
& \quad + \left (\kappa \|\partial_{t} \sigma_{\infty}\|_{L^{2}(0,s;L^{2})} + \lambda_{c} h_{\infty} \|\sigma_{\infty}\|_{L^{2}(0,s;L^{2})} \right ) \|\omega_{k}\|_{L^{2}(0,s;L^{2})} \\
& \leq \frac{D_{0}}{4} \|\nabla \omega_{k}\|_{L^{2}(0,s;L^{2})} + \frac{2D_{1}^{2} \eta^{2}}{D_{0}} \|\nabla \varphi_{k}\|_{L^{2}(0,s;L^{2})}^{2} \\
& \quad + \frac{D_{0}}{4 C_{p}^{2}} \|\omega_{k}\|_{L^{2}(0,s;L^{2})}^{2} + C \left ( 1 + \kappa^{2} \|\partial_{t}\sigma_{\infty}\|_{L^{2}(0,T;L^{2})}^{2} \right ).
\end{aligned}
\end{equation}
Then, substituting \eqref{discrete:apriori:RHS:line4}, \eqref{discrete:apriori:RHS:line1}, \eqref{discrete:apriori:RHS:line2} and \eqref{discrete:apriori:RHS:line3} into \eqref{discrete:apriori:integrated} leads to
\begin{align*}
 \int_{\Omega} & \left [ \frac{\gamma}{\varepsilon} \Psi(\varphi_{k}) + \frac{\gamma \varepsilon}{4} |\nabla \varphi_{k}(s)|^{2} + \mathcal{X} \frac{\kappa}{2} |\omega_{k}|^{2} \right ](s) \, dx + \frac{1}{2} \|\nabla \mu_{k}\|_{L^{2}(0,s;L^{2})}^{2} \\
&  \quad +  \left ( \mathcal{X}\frac{D_{0}}{2} -  2 \chi^{2}(1+C_{p}^{2}) - 2 C_{p}^{2} h_{\infty}^{2} \lambda_{p}^{2} (4 C_{p}^{2} + 1) \right ) \|\nabla \omega_{k}\|_{L^{2}(0,s;L^{2})}^{2} \\
& \leq \left ( C_{p}^{2} + \mathcal{X} \frac{2D_{1}^{2} \eta^{2}}{D_{0}} \right ) \int_{0}^{s} \|\nabla \varphi_{k}\|_{L^{2}}^{2} \, dt \\
& \quad + C \left (1 + \kappa \| \sigma_{0} - \sigma_{\infty}(0)\|_{L^{2}}^{2} + \kappa^{2} \|\partial_{t} \sigma_{\infty}\|_{L^{2}(0,T;L^{2})}^{2} \right ).
\end{align*}
Choosing
\begin{align*}
\mathcal{X} = \frac{4}{D_{0}} \left ( 2 \chi^{2}(1+C_{p}^{2}) + 2 C_{p}^{2} h_{\infty}^{2} \lambda_{p}^{2} (4 C_{p}^{2} + 1) \right ),
\end{align*}
and applying the Gronwall inequality \eqref{Gronwall} with
\begin{align*}
u(s) & := \frac{\gamma}{\varepsilon} \| \Psi(\varphi_{k}) \|_{L^{1}} + \frac{\gamma \varepsilon}{4} \|\nabla \varphi_{k}(s)\|_{L^{2}}^{2}  + \mathcal{X} \frac{\kappa}{2} \|\omega_{k}\|_{L^{2}}^{2}, \\
v(t) & := \frac{1}{2} \|\nabla \mu_{k}\|_{L^{2}}^{2} + \mathcal{X} \frac{D_{0}}{4} \|\nabla \omega_{k}\|_{L^{2}}^{2}, \\
\alpha(s) = \alpha & :=  C \left ( 1 +  \kappa \| \sigma_{0} - \sigma_{\infty}(0)\|_{L^{2}}^{2} +  \kappa^{2} \|\partial_{t} \sigma_{\infty}\|_{L^{2}(0,T;L^{2})}^{2} \right ) \quad \forall s \in (0,T],\\
\beta(t) = \beta & :=  \frac{4}{ \gamma \varepsilon} \left ( C_{p}^{2} + \mathcal{X} \frac{2D_{1}^{2} \eta^{2}}{D_{0}} \right ) \quad \forall t \in (0,T],
\end{align*}
leads to
\begin{align*}
& \sup_{s \in (0,T]} \left ( \|\Psi(\varphi_{k}(s))\|_{L^{1}} + \|\nabla \varphi_{k}(s)\|_{L^{2}}^{2} + \kappa \|\omega_{k}(s)\|_{L^{2}}^{2} \right ) \\
& \qquad + \|\nabla \mu_{k}\|_{L^{2}(0,T;L^{2})}^{2} + \|\nabla \omega_{k}\|_{L^{2}(0,T;L^{2})}^{2} \\
& \quad \leq C \left (1 +  \kappa \| \sigma_{0} - \sigma_{\infty}(0)\|_{L^{2}}^{2} + \kappa^{2} \|\partial_{t} \sigma_{\infty}\|_{L^{2}(0,T;L^{2})}^{2} \right ),
\end{align*}
for some positive constant $C$ independent of $k$ and $\kappa$.  Then, by the Poincar\'{e} inequality, \eqref{sigmakomegakrelation} and \eqref{assump:boundarydata}, we find that
\begin{equation}\label{apriori:first:estimate}
\begin{aligned}
& \sup_{s \in (0,T]} \left ( \|\Psi(\varphi_{k}(s))\|_{L^{1}} + \|\varphi_{k}(s)\|_{H^{1}}^{2} + \kappa \|\sigma_{k}(s)\|_{L^{2}}^{2} \right ) \\
& \qquad +  \|\mu_{k}\|_{L^{2}(0,T;H^{1})}^{2} + \|\sigma_{k}\|_{L^{2}(0,T;H^{1})}^{2} \\
& \quad \leq  C \left (1 +\kappa \| \sigma_{0} - \sigma_{\infty}(0)\|_{L^{2}}^{2} + \kappa^{2} \|\partial_{t} \sigma_{\infty}\|_{L^{2}(0,T;L^{2})}^{2} + \kappa \|\sigma_{\infty}\|_{L^{\infty}(0,T;L^{2})}^{2} \right ),
\end{aligned}
\end{equation}
for some positive constant $C$ independent of $k$ and $\kappa$.  This a priori estimate in turn guarantees that the solution $\{\varphi_{k}, \mu_{k}, \sigma_{k}\}$ to \eqref{Truncated:ODE} exists on the interval $[0,T]$, and thus $t_{k} = T$ for each $k \in \mathbb{N}$.

We now provide some a priori estimates on the time derivatives.  From \eqref{discrete:varphi} and \eqref{discrete:sigma}, we have that
\begin{align*}
\| \partial_{t} \varphi_{k} \|_{L^{2}(0,T;H^{-1})} & \leq \| \nabla \mu_{k} \|_{L^{2}(0,T;L^{2})} + C \left ( 1 + \| \sigma_{k} \|_{L^{2}(0,T;L^{2})} \right ), \\
\kappa \| \partial_{t} \sigma_{k}\|_{L^{2}(0,T;H^{-1})} & \leq C \left ( \| \sigma_{k}\|_{L^{2}(0,T;H^{1})} + \| \nabla \varphi_{k} \|_{L^{2}(0,T;L^{2})} \right ),
\end{align*}
and so by \eqref{apriori:first:estimate},
\begin{equation}
\label{apriori:time}
\begin{aligned}
& \| \partial_{t} \varphi_{k}\|_{L^{2}(0,T;H^{-1})} + \kappa \| \partial_{t} \sigma_{k} \|_{L^{2}(0,T;H^{-1})} \\
& \quad \leq  C \left (1 +\kappa \| \sigma_{0} - \sigma_{\infty}(0)\|_{L^{2}}^{2} + \kappa^{2} \|\partial_{t} \sigma_{\infty}\|_{L^{2}(0,T;L^{2})}^{2} + \kappa \|\sigma_{\infty}\|_{L^{\infty}(0,T;L^{2})}^{2} \right ),
\end{aligned}
\end{equation}
for some positive constant $C$ independent of $k$ and $\kappa$.

\subsubsection{Passing to the limit}\label{sec:passinglimit}
From \eqref{apriori:first:estimate} and \eqref{apriori:time} we find that
\begin{align*}
\{\varphi_{k}\}_{k \in \mathbb{N}} & \text{ is bounded in } L^{\infty}(0,T;H^{1}) \cap H^{1}(0,T;H^{-1}), \\
\{\sigma_{k}\}_{k \in \mathbb{N}} & \text{ is bounded in } L^{\infty}(0,T;L^{2}) \cap L^{2}(0,T;H^{1}) \cap H^{1}(0,T;H^{-1}), \\
\{\mu_{k}\}_{k \in \mathbb{N} } & \text{ is bounded in } L^{2}(0,T;H^{1}),
\end{align*}
and by the standard compactness results we obtain for a non-relabeled subsequence the following convergences
\begin{alignat*}{3}
\varphi_{k} & \rightarrow \varphi && \quad \text{ weakly-}* && \quad \text{ in } L^{\infty}(0,T;H^{1}) \cap H^{1}(0,T;H^{-1}), \\
\varphi_{k} & \rightarrow \varphi && \quad \text{ strongly } && \quad \text{ in } C^{0}([0,T];H^{r}) \cap L^{2}(0,T;L^{m}) \text{ and a.e. in } \Omega \times (0,T), \\
\sigma_{k} & \rightarrow \sigma && \quad \text{ weakly-}* && \quad \text{ in } L^{2}(0,T;H^{1}) \cap L^{\infty}(0,T;L^{2}) \cap H^{1}(0,T;H^{-1}), \\
\sigma_{k} & \rightarrow \sigma && \quad \text{ strongly } && \quad \text{ in } L^{2}(0,T;L^{m}) \text{ and a.e. in } \Omega \times (0,T), \\
\mu_{k} & \rightarrow \mu && \quad \text{ weakly } && \quad \text{ in } L^{2}(0,T;H^{1}),
\end{alignat*}
for $0 \leq r < 1$ and $m \in [1,6)$.

Fix $j$ and consider $\delta(t) \in C^{\infty}_{c}(0,T)$.  Then, $\delta(t)w_{j} \in L^{2}(0,T;H^{2} \cap H^{1}_{0})$.  We multiple \eqref{Truncated:discrete} with $\delta(t)$ and integrate in time from $0$ to $T$, leading to
\begin{subequations}\label{Truncated:discrete:time}
\begin{align}
0 & = \int_{0}^{T} \int_{\Omega} \delta(t) \left ( \partial_{t}\varphi_{k} w_{j} + \nabla \mu_{k} \cdot \nabla w_{j} - (\lambda_{p} \sigma_{k} - \lambda_{a}) h(\varphi_{k}) w_{j} \right )\, dx \, dt, \label{discrete:varphi:time} \\
0 & = \int_{0}^{T} \int_{\Omega}  \delta(t) \left ( \mu_{k} w_{j} - \frac{\gamma}{\varepsilon} \Psi'(\varphi_{k}) w_{j} - \gamma \varepsilon \nabla \varphi_{k} \cdot \nabla w_{j} + \chi \sigma_{k} w_{j} \right ) \, dx \, dt, \label{discrete:mu:time} \\
0 & = \int_{0}^{T} \int_{\Omega}  \delta(t) \left ( \kappa \partial_{t} \sigma_{k} w_{j} + D(\varphi_{k}) (\nabla \sigma_{k} - \eta \nabla \varphi_{k}) \cdot \nabla w_{j} + \lambda_{c} \sigma_{k} h(\varphi_{k}) w_{j} \right ) \,dx \, dt. \label{discrete:sigma:time}
\end{align}
\end{subequations}

By continuity of $h(\cdot)$, we see that $h(\varphi_{k}) \to h(\varphi)$ a.e. in $\Omega \times (0,T)$.  Thanks to the fact that $h(\cdot)$ is bounded, applying Lebesgue dominated convergence theorem to $(h(\varphi_{k}) - h(\varphi)) \delta w_{j}$ yields
\begin{align*}
\|(h(\varphi_{k}) - h(\varphi)) \delta w_{j}\|_{L^{2}(0,T;L^{2})} \to 0 \text{ as } k \to \infty.
\end{align*}
Together with the weak convergence $\sigma_{k} \rightharpoonup \sigma$ in $L^{2}(0,T;L^{2})$, we obtain by the product of weak-strong convergence
\begin{align*}
\int_{0}^{T} \int_{\Omega} h(\varphi_{k}) \sigma_{k} \delta w_{j} \, dx \, dt \to \int_{0}^{T} \int_{\Omega} h(\varphi) \sigma \delta w_{j} \, dx \, dt \text{ as } k \to \infty.
\end{align*}
The terms involving $D(\cdot)$ can be treated in a similar fashion.  The remaining non-trivial part is to pass to the limit in the term involving $\Psi'$.  By \eqref{assump:Psi} and \eqref{apriori:first:estimate}, we have that $\{\Psi'(\varphi_{k})\}_{k \in \mathbb{N}}$ is bounded uniformly in $L^{\infty}(0,T;L^{s})$ for some $s \in (1,2]$.  Thus, there exists a function $\xi \in L^{\infty}(0,T;L^{s})$ such that
\begin{align*}
\Psi'(\varphi_{k}) \to \xi \text{ weakly-}* \text{ in } L^{\infty}(0,T;L^{s}).
\end{align*}
By the continuity of $\Psi'$ and the a.e. convergence of $\varphi_{k}$ to $\varphi$, we have $\Psi'(\varphi_{k}) \to \Psi'(\varphi)$ a.e. in $\Omega \times (0,T)$.  The uniqueness of a.e. and weak limits shows that $\xi = \Psi'(\varphi)$.  We now pass to the limit $k \to \infty$ in \eqref{Truncated:discrete:time}, and the application of the aforementioned convergence results leads to
\begin{align*}
0 & = \int_{0}^{T} \delta(t) \left ( \langle \partial_{t}\varphi, w_{j} \rangle  +  \int_{\Omega} \nabla \mu \cdot \nabla w_{j} - (\lambda_{p} \sigma - \lambda_{a}) h(\varphi) w_{j} \, dx \right ) \, dt,  \\
0 & = \int_{0}^{T}   \delta(t) \left ( \int_{\Omega} \mu w_{j} - \frac{\gamma}{\varepsilon} \Psi'(\varphi) w_{j} - \gamma \varepsilon \nabla \varphi \cdot \nabla w_{j} + \chi \sigma w_{j} \,dx \right ) \,dt,  \\
0 & = \int_{0}^{T}  \delta(t) \left ( \langle \kappa \partial_{t} \sigma, w_{j} \rangle +  \int_{\Omega} D(\varphi) (\nabla \sigma - \eta \nabla \varphi) \cdot \nabla w_{j} + \lambda_{c} \sigma h(\varphi) w_{j} \, dx \right ) \, dt,
\end{align*}
which holds for arbitrary $\delta(t) \in C^{\infty}_{c}(0,T)$.  Hence, we have for a.e. $t \in (0,T)$,
\begin{subequations}
\begin{align}
\langle \partial_{t}\varphi, w_{j} \rangle & = \int_{\Omega} - \nabla \mu \cdot \nabla w_{j} + (\lambda_{p} \sigma - \lambda_{a}) h(\varphi) w_{j} \, dx, \\
\int_{\Omega} \mu w_{j} \, dx & = \int_{\Omega} \frac{\gamma}{\varepsilon} \Psi'(\varphi) w_{j} + \gamma \varepsilon \nabla \varphi \cdot \nabla w_{j} - \chi \sigma w_{j} \, dx, \label{limit:mu} \\
\langle \kappa \partial_{t} \sigma, w_{j}\rangle & = -\int_{\Omega} D(\varphi) (\nabla \sigma - \eta \nabla \varphi) \cdot \nabla w_{j} + \lambda_{c} \sigma h(\varphi) w_{j} \, dx. \label{limit:sigma}
\end{align}
\end{subequations}
This holds for all $j \geq 1$, and as $\{ w_{j}\}_{j \in \mathbb{N}}$ is a basis for $H^{1}_{0}$, we see that $\{ \varphi, \mu, \sigma \}$ satisfy \eqref{truncated:weak} for all $\zeta$, $\lambda$, $\xi \in H^{1}_{0}$.  Moreover, from the strong convergence of $\varphi_{k}$ to $\varphi$ in $C^{0}([0,T];L^{2})$, and the fact that $\varphi_{k}(0) \to \varphi_{0}$ in $L^{2}$ it holds that $\varphi(0) = \varphi_{0}$. Meanwhile, by the continuous embedding
\begin{align*}
L^{2}(0,T;H^{1}) \cap H^{1}(0,T;H^{-1}) \subset C^{0}([0,T];L^{2}),
\end{align*}
and that $\sigma_{k}(0) \to \sigma_{0}$ in $L^{2}$, we have $\sigma(0) = \sigma_{0}$.  This shows that $\{\varphi, \mu, \sigma\}$ is a weak solution of \eqref{Intro:CHN}.

\begin{remark}
For the reader's convenience, we outline an alternative argument for showing
\begin{align}\label{Psi':conv:p<5}
\int_{0}^{T} \int_{\Omega} \Psi'(\varphi_{k}) \delta w_{j} \, dx \, dt \to \int_{0}^{T} \int_{\Omega} \Psi'(\varphi) \delta w_{j} \, dx \, dt \text{ as } k \to \infty
\end{align}
when the potential $\Psi$ satisfies \eqref{Psi'5thorder} instead of \eqref{assump:Psi} that uses the generalized Lebesgue dominated convergence theorem.  For any $q < 6$, choose $r = \frac{3}{2} - \frac{3}{q} < 1$.  Then, the Sobolev embedding $H^{r} \subset L^{q}$ yields that $C^{0}([0,T];H^{r}) \subset L^{q}(0,T;L^{q})$ for $q < 6$, and by the above compactness results, we have
$\varphi_{k} \to \varphi$ strongly in $L^{q}(0,T;L^{q}) \cong L^{q}(\Omega \times (0,T))$ for $q < 6$.  For $p \in [1,5)$, let $q = \frac{6}{5}p < 6$.  Then, by H\"{o}lder's inequality,
\begin{align*}
\int_{0}^{T} \int_{\Omega} | |\varphi_{k} - \varphi|^{p} \delta w_{j}| \, dx \, dt & \leq \|\varphi_{k} - \varphi\|_{L^{q}(0,T;L^{q})}^{p} \|w_{j}\|_{L^{6}} \|\delta\|_{L^{6}(0,T)} \to 0 \text{ as } k \to \infty.
\end{align*}
Hence, $|\varphi_{k} - \varphi|^{p} \delta w_{j} \to 0$ strongly in $L^{1}(\Omega \times (0,T))$.  A short computation shows that
\begin{itemize}
\item $|\varphi_{k}|^{p} |\delta w_{j}| \leq C(p) (|\varphi|^{p} + |\varphi_{k} - \varphi|^{p}) |\delta w_{j}|  \in L^{1}(\Omega \times (0,T))$ for all $k$,
\item $(|\varphi|^{p} + |\varphi_{k} - \varphi|^{p}) |\delta w_{j}| \to |\varphi|^{p} |\delta w_{j}|$ a.e in  $\Omega \times (0,T)$ as $k \to \infty$,
\item $ \int_{0}^{T} \int_{\Omega} (|\varphi|^{p} + |\varphi_{k} - \varphi|^{p}) |\delta w_{j}| \, dx \, dt \to \int_{0}^{T} \int_{\Omega} |\varphi|^{p} |\delta w_{j}| \, dx \, dt$ as $k \to \infty$.
\end{itemize}
By the generalized Lebesgue dominated convergence theorem, we have
\begin{align*}
\int_{0}^{T} \int_{\Omega} (1 +  |\varphi_{k}|^{p}) \delta w_{j} \, dx \, dt \to \int_{0}^{T} \int_{\Omega} ( 1 + |\varphi|^{p}) \delta w_{j} \, dx \, dt \text{ as } k \to \infty \text{ for } p \in [1,5).
\end{align*}
Using that $\Psi'(\varphi_{k}) \to \Psi'(\varphi)$ a.e. in $\Omega \times (0,T)$ as $k \to \infty$, and
\begin{align*}
|\Psi'(\varphi_{k}) \delta w_{j}| \leq k_{2}(1 + |\varphi_{k}|^{p} ) |\delta w_{j}| \in L^{1}(\Omega \times (0,T)) \; \forall k \in \mathbb{N} \text{ and } p \in [1,5),
\end{align*}
and the generalized Lebesgue dominated convergence theorem again, we obtain \eqref{Psi':conv:p<5}.
\end{remark}

\subsection{Energy inequality}\label{sec:energyineq}
By the a.e. convergence $\varphi_{k} \to \varphi$ in $\Omega \times (0,T)$, the non-negativity and continuity of $\Psi$, Fatou's lemma and the estimate \eqref{apriori:first:estimate}, we have that for a.e. $s \in (0,T]$,
\begin{align*}
&\int_{\Omega} \Psi(\varphi(s)) \, dx \leq \liminf_{k \to \infty} \int_{\Omega} \Psi(\varphi_{k}(s)) \, dx \\
& \leq  C \left (1 + \kappa \|\sigma_{0} - \sigma_{\infty}(0)\|_{L^{2}}^{2} + \kappa^{2} \| \partial_{t} \sigma_{\infty}\|_{L^{2}(0,T;L^{2})}^{2} + \kappa \|\sigma_{\infty}\|_{L^{\infty}(0,T;L^{2})}^{2} \right ).
\end{align*}
Taking the supremum in $s \in (0,T]$ in the above, and using the following convergences:
\begin{alignat*}{3}
\varphi_{k} & \rightarrow \varphi && \quad \text{ weakly-}* && \quad \text{ in } L^{\infty}(0,T;H^{1}), \\
\sigma_{k} & \rightarrow \sigma && \quad \text{ weakly-}* && \quad \text{ in } L^{2}(0,T;H^{1}) \cap L^{\infty}(0,T;L^{2}), \\
\mu_{k} & \rightarrow \mu && \quad \text{ weakly } && \quad \text{ in } L^{2}(0,T;H^{1}),
\end{alignat*}
along with the weak/weak-$*$ lower semicontinuity of the Sobolev norms, we find that passing to the limit $k \to \infty$ in \eqref{apriori:first:estimate} leads to \eqref{energyineq}.

\subsection{Further regularity}\label{sec:furtherreg}
Note that \eqref{truncated:weak:mu} is the weak formulation of
\begin{subequations}\label{elliptic:reg:varphi}
\begin{alignat}{2}
- \gamma \varepsilon \Delta \varphi  & = \mu - \frac{\gamma}{\varepsilon} \Psi'(\varphi) + \chi \sigma  && \text{ in } \Omega, \label{elliptic:reg:varphi:equ} \\
 \varphi & = -1 && \text{ on } \Gamma.
\end{alignat}
\end{subequations}
To prove the $L^{2}(0,T;W^{2,6})$-regularity assertion, it suffices to show that $\Psi'(\varphi)$ is bounded in $L^{6}$, as $\mu + \chi \sigma$ already belongs to $H^{1}\subset L^{6}$.  We define iteratively, a sequence of numbers $\{l_{j}\}_{j \in \mathbb{N}}$, by
\begin{align}\label{Regularity:exponents}
l_{1} \geq 1, \quad
l_{1}p \leq 6, \quad l_{j+1} = \frac{6 l_{j}}{6  - (5-p) l_{j} },
\end{align}
where $p \in (1,5)$ is the exponent in \eqref{Psi'5thorder}.  While $l_{j}$ and $l_{j+1}$ are positive and greater than $1$, by the Gagliardo--Nirenberg inequality, it holds that
\begin{align*}
\|f\|_{L^{p l_{j+1}}}^{2p} \leq C \|f\|_{W^{2,l_{j}}}^{2} \|f\|_{L^{6}}^{2p-2} \Leftrightarrow \frac{1}{l_{j+1}p} =\left ( \frac{1}{l_{j}} - \frac{2}{3} \right )\frac{1}{p} + \frac{1}{6} - \frac{1}{6p}.
\end{align*}
Rearranging the equality on the right-hand side yields exactly the relation $ l_{j+1} = \frac{6 l_{j}}{6  - (5-p) l_{j} }$, and so
\begin{align}\label{Regularity:GN:Bochner}
L^{2}(0,T;W^{2,l_{j}}) \cap L^{\infty}(0,T;L^{6}) \subset L^{2p}(0,T;L^{l_{j+1}p}).
\end{align}
The bootstrapping procedure is as follows: 

\vspace*{4pt}\noindent\textbf{First step.} By \eqref{Psi'5thorder} and \eqref{Regularity:exponents}, we have
\begin{align*}
\|\Psi'(\varphi)\|_{L^{l_{1}}}^{2} \leq C(k_{1}) \left (1+\|\varphi\|_{L^{6}}^{2 p} \right ) \leq C(k_{1}, C_{s}, l_{1}, p) \left (1 + \|\varphi\|_{H^{1}}^{2 p} \right ),
\end{align*}
and due to the boundedness of $\varphi$ in $L^{\infty}(0,T;H^{1})$, the right-hand side of \eqref{elliptic:reg:varphi:equ} belongs to $L^{2}(0,T;L^{l_{1}})$ and elliptic regularity theory yields that $\varphi \in L^{2}(0,T;W^{2,l_{1}})$. 

\vspace*{4pt}\noindent\textbf{$j$-th step.} For each $j > 1$, from \eqref{Psi'5thorder} and \eqref{Regularity:GN:Bochner} we have that
\begin{align*}
\|\Psi'(\varphi)\|_{L^{l_{j}}}^{2} \leq C(k_{1}) \left ( 1 + \|\varphi\|_{L^{l_{j}p}}^{2p} \right ) \leq C \left ( 1 + \|\varphi\|_{W^{2,l_{j-1}}}^{2} \|\varphi\|_{L^{6}}^{2p-2} \right ),
\end{align*}
which in turn implies that $\Psi'(\varphi)$ belongs to $L^{2}(0,T;L^{l_{j}})$, and the application of elliptic regularity then gives that $\varphi \in L^{2}(0,T;W^{2,l_{j}})$.

\vspace*{4pt}

Let us point out that the function $g(x) := \frac{6x}{6-(5-p)x}$ is strictly increasing and positive in the interval $[1,\frac{6}{5-p})$.  In particular,
\begin{align*}
l_{j} < l_{j+1} & \text{ if } l_{j}, l_{j+1} < \frac{6}{5-p}, \\
l_{j+1} < 0 & \text{ if } l_{j} > \frac{6}{5-p}.
\end{align*}
By the strictly increasing property, after a finite number of steps, we have $l_{j} \geq \frac{6}{6-p}$, and at this point the bootstrapping procedure is terminated.  Then, $l_{j+1} \geq 6$ and \eqref{Psi'5thorder} and \eqref{Regularity:GN:Bochner} imply that $\varphi$ is bounded in $L^{2p}(0,T;L^{6p})$ and $\Psi'(\varphi)$ is bounded in $L^{2}(0,T;L^{6})$.  For example, in the case $p = 3$ which corresponds to the double-well potential $(s^{2}-1)^{2}$, we have $l_{1} = 2$ and $l_{2} = 6$, and so the bootstrapping procedure ends in two steps.  In the case $p = 4$, we have $l_{1} = \frac{3}{2}$, $l_{2} = 2$, $l_{3} = 3$ and $l_{4} = 6$, and the bootstrapping procedure ends in four steps.

Now let $\Gamma$ be of class $C^{3}$, we show that $\Psi'(\varphi) \in L^{2}(0,T;H^{1})$.  Then, the right-hand side of \eqref{elliptic:reg:varphi:equ} belongs to $H^{1}$ for a.e. $t \in (0,T)$ and by elliptic regularity theory, one obtains that $\varphi \in L^{2}(0,T;H^{3})$.  From the growth assumptions \eqref{assump:Psi''} on $\Psi'$ and $\Psi''$, we find that
\begin{align*}
|\Psi'(\varphi)|^{2} \leq 2k_{2}^{2} (1 + |\varphi|^{2p}), \quad |\Psi''(\varphi) \nabla \varphi|^{2} \leq 2k_{3}^{2}( 1 + |\varphi|^{2p-2}) |\nabla \varphi|^{2}
\end{align*}
for $p < 5$.  In dimension $d = 3$, the Gagliardo--Nirenberg inequality yields that (substituting $l_{j+1} = 2$ in \eqref{Regularity:exponents} and rearrange the equation for $l_{j}$)
\begin{align*}
\|\varphi\|_{L^{2p}}^{2p} \leq C \|\varphi\|_{W^{2,m}}^{2} \|\varphi\|_{L^{6}}^{2p-2} \text{ for } m = \frac{6}{8 - p}.
\end{align*}
Then, for any $p < 5$, the $L^{2}(0,T;W^{2,6})$-regularity established above yields that
\begin{align*}
L^{2}(0,T;W^{2,6}) \cap L^{\infty}(0,T;H^{1}) \subset L^{2p}(0,T;L^{2p})  \Rightarrow \Psi'(\varphi) \in L^{2}(0,T;L^{2}).
\end{align*}
Similarly, the Gagliardo--Nirenberg inequality yields that for $p < 5$,
\begin{align*}
\|\varphi\|_{L^{\infty}}^{2p-2} \leq C \|\varphi\|_{W^{2,q}}^{2} \|\varphi\|_{L^{6}}^{2p-4} \text{ for } q = \frac{6}{6-p}.
\end{align*}
From the $L^{2}(0,T;W^{2,6})$-regularity, we have that $\varphi \in L^{2p-2}(0,T;L^{\infty})$ and thus
\begin{align*}
\|\Psi''(\varphi) \nabla \varphi\|_{L^{2}(0,T;L^{2})}^{2} \leq C \left ( 1 + \|\varphi\|_{L^{2p-2}(0,T;L^{\infty})}^{2p-2}\right ) \|\nabla \varphi\|_{L^{\infty}(0,T;L^{2})}^{2} ,
\end{align*}
and this establishes that $\Psi'(\varphi) \in L^{2}(0,T;H^{1})$.
\begin{remark}
The restriction to $L^{2}(0,T;W^{2,6})$ for the first regularity assertion of Theorem \ref{thm:truncated:regularity} is due to the fact that $H^{3}$-regularity requires a $C^{3}$-boundary.  If $\Gamma$ is only a $C^{2}$-boundary, then at best we can only deduce $\varphi \in L^{2}(0,T;W^{2,6})$ even if the right-hand side of \eqref{elliptic:reg:varphi:equ} belongs to $L^{2}(0,T;H^{1})$.
\end{remark}

\subsection{Continuous dependence}\label{sec:ctsdep}
Suppose $D(\cdot) = D_{0}$ is constant, $h(\cdot)$ is Lipschitz continuous with Lipschitz constant $\mathrm{L}_{h}$, and $\Psi'$ satisfies \eqref{Psi'':uniqueness}.  Given two regular weak solution triplets $\{ \varphi_{i}, \mu_{i}, \sigma_{i}\}_{i = 1, 2}$ to \eqref{Intro:CHN} with corresponding initial data $\{\varphi_{i,0}, \sigma_{i,0}\}_{i = 1,2}$ and boundary data $\{\mu_{\infty,i}, \sigma_{\infty,i} \}_{i = 1,2}$, let $\varphi := \varphi_{1} - \varphi_{2}$, $\mu:= \mu_{1} - \mu_{2}$, $\sigma := \sigma_{1} - \sigma_{2}$, $\mu_{\infty} := \mu_{\infty,1} - \mu_{\infty,2}$, and $\sigma_{\infty} := \sigma_{\infty,1} - \sigma_{\infty,2}$ denote the differences, respectively. Then, we have that
\begin{align*}
\varphi & \in L^{\infty}(0,T;H^{1}_{0}) \cap L^{2}(0,T;H^{3}) \cap H^{1}(0,T;H^{-1}), \\
\mu & \in \mu_{\infty} + L^{2}(0,T;H^{1}_{0}), \\
\sigma & \in \left ( \sigma_{\infty} + L^{2}(0,T;H^{1}_{0}) \right ) \cap L^{\infty}(0,T;L^{2}) \cap H^{1}(0,T;H^{-1}),
\end{align*}
and they satisfy $\varphi(0) = \varphi_{0} := \varphi_{1,0} - \varphi_{2,0}$, $\sigma(0) = \sigma_{0} := \sigma_{1,0} - \sigma_{2,0}$, and
\begin{align*}
0 & = \langle \partial_{t}\varphi, \zeta \rangle + \int_{\Omega} \nabla \mu \cdot \nabla \zeta + ( ( \lambda_{a} - \lambda_{p} \sigma_{1}  ) (h(\varphi_{1}) - h(\varphi_{2})) - \lambda_{p} \sigma h(\varphi_{2}))  \zeta \, dx, \\
0 & = \int_{\Omega} \mu \lambda - \frac{\gamma}{\varepsilon} (\Psi'(\varphi_{1}) - \Psi'(\varphi_{2})) \lambda - \gamma \varepsilon \nabla \varphi \cdot \nabla \lambda + \chi \sigma \lambda \, dx, \\
0 & = \langle \kappa \partial_{t} \sigma, \xi \rangle + \int_{\Omega} D_{0} (\nabla \sigma - \eta \nabla \varphi) \cdot \nabla \xi + \lambda_{c} (\sigma_{1} (h(\varphi_{1}) - h(\varphi_{2})) + \sigma h(\varphi_{2})) \xi \, dx,
\end{align*}
for all $\zeta, \lambda, \xi \in H^{1}_{0}(\Omega)$ and a.e. $t \in (0,T)$.  In the following, we will often use
\begin{align}\label{Poincare:sigmaL2}
\|\sigma\|_{L^{2}} \leq \|\sigma - \sigma_{\infty}\|_{L^{2}} + \|\sigma_{\infty}\|_{L^{2}} \leq C_{p} \|\nabla \sigma\|_{L^{2}} + C(C_{p}) \|\sigma_{\infty}\|_{H^{1}},
\end{align}
obtained from applying the Poincar\'{e} inequality to $\sigma - \sigma_{\infty}$.  For positive constant $\mathcal{Z}$ yet to be determined, substituting $\xi = \mathcal{Z} (\sigma - \sigma_{\infty})$ leads to
\begin{equation}\label{uniqueness:sigmaequ:testedsigma}
\begin{aligned}
\mathcal{Z} \frac{\kappa}{2} \frac{d}{dt} \|\sigma\|_{L^{2}}^{2} & + \mathcal{Z} D_{0} \|\nabla \sigma\|_{L^{2}}^{2} + \mathcal{Z} \lambda_{c} \int_{\Omega} \underbrace{h(\varphi_{2}) |\sigma|^{2}}_{\geq 0} \, dx \\
= & \, \mathcal{Z} \kappa \langle \partial_{t}\sigma, \sigma_{\infty} \rangle + \mathcal{Z} \int_{\Omega}  D_{0} \eta \nabla \varphi \cdot \nabla (\sigma - \sigma_{\infty}) + D_{0} \nabla \sigma \cdot \nabla \sigma_{\infty} \, dx \\
&  - \mathcal{Z} \int_{\Omega} \lambda_{c} \sigma_{1} (h(\varphi_{1}) - h(\varphi_{2})) (\sigma - \sigma_{\infty}) - \lambda_{c} h(\varphi_{2}) \sigma \sigma_{\infty} \, dx .
\end{aligned}
\end{equation}
Then, substituting $\zeta = \gamma \varepsilon \varphi$, $\lambda = \mu - \mu_{\infty}$ and $\lambda = \mathcal{Y} \varphi$, for some positive constant $\mathcal{Y}$ yet to be determined, and summing with \eqref{uniqueness:sigmaequ:testedsigma}, integrating over $t$ from $0$ to $s \in (0,T]$ leads to
\begin{equation}\label{CHN:uniqueness:apriori:est}
\begin{aligned}
& \left ( \frac{\gamma \varepsilon}{2} \|\varphi(s)\|_{L^{2}}^{2} + \mathcal{Z} \frac{\kappa}{2} \|\sigma(s)\|_{L^{2}}^{2} \right ) + \int_{0}^{s} \|\mu\|_{L^{2}}^{2} + \mathcal{Z} D_{0} \|\nabla \sigma\|_{L^{2}}^{2} + \mathcal{Y} \gamma \varepsilon \|\nabla \varphi\|_{L^{2}}^{2} \, dt \\
& \quad \leq J_{1} + J_{2} + J_{3} + J_{4} + J_{5} +  \left ( \frac{\gamma \varepsilon}{2} \|\varphi_{0}\|_{L^{2}}^{2} +  \mathcal{Z} \frac{\kappa}{2} \|\sigma_{0}\|_{L^{2}}^{2} \right ),
\end{aligned}
\end{equation}
where
\begin{align*}
&J_{1}  = - \mathcal{Z} \kappa \int_{0}^{s} \int_{\Omega} \sigma \partial_{t} \sigma_{\infty} \, dx \, dt + \mathcal{Z} \kappa \int_{\Omega} \sigma(s)  \sigma_{\infty}(s) - \sigma_{0} \sigma_{\infty}(0) \, dx, \\
&J_{2} \\& = \int_{0}^{s} \int_{\Omega} \frac{\gamma}{\varepsilon} (\Psi'(\varphi_{1}) - \Psi'(\varphi_{2}))(\mu - \mu_{\infty} - \mathcal{Y} \varphi ) + \mathcal{Y}(\mu + \chi \sigma) \varphi - \gamma \varepsilon \nabla \varphi \cdot \nabla \mu_{\infty} \, dx \, dt, \\
&J_{3}  = \int_{0}^{s} \int_{\Omega} \mathcal{Z} D_{0} \eta \nabla \varphi \cdot \nabla (\sigma- \sigma_{\infty}) + \mathcal{Z} D_{0} \nabla \sigma \cdot \nabla \sigma_{\infty} - \chi \sigma (\mu - \mu_{\infty}) + \mu \mu_{\infty} \, dx \, dt, \\
&J_{4}  = \int_{0}^{s} \int_{\Omega} -\gamma \varepsilon \lambda_{a} (h(\varphi_{1}) - h(\varphi_{2}))\varphi +  \lambda_{p} \gamma \varepsilon h(\varphi_{2}) \sigma \varphi + \mathcal{Z} \lambda_{c} h(\varphi_{2}) \sigma \sigma_{\infty} \, dx  \, dt, \\
&J_{5}  = \int_{0}^{s} \int_{\Omega}  (\gamma \varepsilon \lambda_{p} \varphi + \lambda_{c} \mathcal{Z} (\sigma_{\infty} -  \sigma)) \sigma_{1} (h(\varphi_{1}) - h(\varphi_{2})) \, dx \, dt.
\end{align*}

Application of H\"{o}lder's inequality and Young's inequality leads to
\begin{subequations}\label{ctsdep:J1J3J4}
\begin{align}
\label{ctsdep:J1} |J_{1}| & \leq  \frac{1}{4} \|\sigma\|_{L^{2}(0,s;L^{2})}^{2} + \kappa^{2} \mathcal{Z}^{2} \| \partial_{t} \sigma_{\infty}\|_{L^{2}(0,T;L^{2})}^{2} + \mathcal{Z} \frac{\kappa}{4} \|\sigma(s)\|_{L^{2}}^{2} \\
\notag & \quad + 2 \kappa \mathcal{Z} \|\sigma_{\infty}\|_{L^{\infty}(0,T;L^{2})}^{2} + \mathcal{Z} \frac{\kappa}{4} \|\sigma_{0}\|_{L^{2}}^{2}, \\
\label{ctsdep:J3} |J_{3}| & \leq \frac{2 D_{0}}{4} \mathcal{Z} \|\nabla \sigma\|_{L^{2}(0,s;L^{2})}^{2} + \frac{2}{8} \|\mu\|_{L^{2}(0,s;L^{2})}^{2} + 2  \mathcal{Z} D_{0} \eta^{2} \|\nabla \varphi\|_{L^{2}(0,s;L^{2})}^{2}  \\
\notag & \quad  + \frac{5 D_{0}}{4} \mathcal{Z} \|\nabla \sigma_{\infty}\|_{L^{2}(0,T;L^{2})}^{2} + (2 \chi^{2} + 1) \|\sigma\|_{L^{2}(0,s;L^{2})}^{2} + \frac{9}{4} \|\mu_{\infty}\|_{L^{2}(0,T;L^{2})}^{2}, \\
\label{ctsdep:J4} |J_{4}| & \leq  \left ( \gamma \varepsilon \lambda_{a} \mathrm{L}_{h} + (\lambda_{p} \gamma \varepsilon h_{\infty})^{2} \right ) \|\varphi\|_{L^{2}(0,s;L^{2})}^{2} + \frac{2}{4} \|\sigma\|_{L^{2}(0,s;L^{2})}^{2} \\
\notag & \quad + (\lambda_{c} h_{\infty} \mathcal{Z})^{2} \|\sigma_{\infty}\|_{L^{2}(0,T;L^{2})}^{2}.
\end{align}
\end{subequations}
By \eqref{Psi'':uniqueness} and the following Gagliardo--Nirenberg interpolation inequality in three-dimensions,
\begin{align*}
\|f\|_{L^{\infty}} \leq C \|f\|_{H^{3}}^{\frac{1}{4}} \|f\|_{L^{6}}^{\frac{3}{4}}  \leq \hat{C} \|f\|_{H^{3}}^{\frac{1}{4}} \|f\|_{H^{1}}^{\frac{3}{4}},
\end{align*}
for positive constants $C, \hat{C}$ depending only on $\Omega$, we find that
\begin{align*}
& \left |\int_{\Omega} (\Psi'(\varphi_{1}) - \Psi'(\varphi_{2}))(\mu - \mu_{\infty} - \mathcal{Y} \varphi) \, dx \right | \\
& \quad \leq k_{3} \int_{\Omega} (1 + |\varphi_{1}|^{4} + |\varphi_{2}|^{4}) (\mathcal{Y} |\varphi|^{2} + |\varphi| |\mu| + |\varphi| |\mu_{\infty}|) \, dx \\
& \quad \leq k_{3} \left  ( 1+ \|\varphi_{1}\|_{L^{\infty}}^{4} + \|\varphi_{2}\|_{L^{\infty}}^{4} \right ) \left (\mathcal{Y} \|\varphi\|_{L^{2}}^{2} + \|\varphi\|_{L^{2}} \left ( \|\mu\|_{L^{2}} + \|\mu_{\infty}\|_{L^{2}} \right ) \right ) \\
& \quad \leq C \big{(} 1+ \textstyle \sum_{i=1,2} \|\varphi_{i}\|_{H^{3}} \|\varphi_{i}\|_{H^{1}}^{3} \big{)} \left ( \mathcal{Y} \|\varphi\|_{L^{2}}^{2} + \|\varphi\|_{L^{2}} \left (\|\mu\|_{L^{2}} + \|\mu_{\infty}\|_{L^{2}} \right ) \right ),
\end{align*}
where $C$ depends only on $k_{3}$ and $\hat{C}$.  Then, we can estimate $J_{2}$ as follows,
\begin{equation*}
\begin{aligned}
|J_{2}| & \leq  C \int_{0}^{s} \left ( 1 + \textstyle \sum_{i=1,2} \|\varphi_{i}\|_{H^{3}} \|\varphi_{i}\|_{H^{1}}^{3} \right ) \mathcal{Y} \|\varphi\|_{L^{2}}^{2} \, dx 
\end{aligned}
\end{equation*}\begin{equation}\label{ctsdep:J2}
\begin{aligned}& \quad + C \int_{0}^{s} \left ( 1 + \textstyle \sum_{i=1,2} \|\varphi_{i}\|_{H^{3}}^{2} \|\varphi_{i}\|_{H^{1}}^{6} \right ) \|\varphi\|_{L^{2}}^{2} \, dt  + \frac{2}{8} \|\mu\|_{L^{2}(0,s;L^{2})}^{2}  \\
& \quad+ \|\mu_{\infty}\|_{L^{2}(0,T;L^{2})}^{2} + \frac{\chi^{2}}{4} \|\sigma\|_{L^{2}(0,s;L^{2})}^{2} + 3 \mathcal{Y}^{2} \|\varphi\|_{L^{2}(0,s;L^{2})}^{2} \\
& \quad  + \gamma \varepsilon \|\nabla \varphi\|_{L^{2}(0,s;L^{2})}^{2}+ \frac{\gamma \varepsilon}{4} \|\nabla \mu_{\infty}\|_{L^{2}(0,T;L^{2})}^{2}.
\end{aligned}
\end{equation}

By the Lipschitz continuity of $h(\cdot)$ and \eqref{Poincare:sigmaL2}, we have
\begin{equation}\label{ctsdep:sigmahvarphi}
\begin{aligned}
|J_{5}| & = \left | \int_{0}^{s} \int_{\Omega}  \sigma_{1} (h(\varphi_{1}) - h(\varphi_{2})) \left ( \mathcal{Z} \lambda_{c} (\sigma_{\infty} - \sigma) + \gamma \varepsilon \lambda_{p} \varphi \right ) \, dx \, dt \right | \\
& \leq \mathrm{L}_{h} \int_{0}^{s} \|\sigma_{1}\|_{L^{2}} \|\varphi\|_{L^{4}} \left (  \mathcal{Z} \lambda_{c} \left ( \|\sigma\|_{L^{4}} + \|\sigma_{\infty}\|_{L^{4}} \right ) + \gamma \varepsilon \lambda_{p} \|\varphi\|_{L^{4}} \right ) \,dt \\
& \leq \overline{C} \|\varphi\|_{L^{2}(0,s;H^{1})} \mathcal{Z} \lambda_{c} \left ( \|\sigma\|_{L^{2}(0,s;H^{1})} + \|\sigma_{\infty}\|_{L^{2}(0,T;H^{1})} \right ) \\
& \quad + \overline{C} \gamma \varepsilon \lambda_{p} \|\varphi\|_{L^{2}(0,s;H^{1})}^{2} \\
& \leq \overline{C} \left ( 2\overline{C} \lambda_{c}^{2} \mathcal{Z}^{2}  + \gamma \varepsilon \lambda_{p} \right ) \|\varphi\|_{L^{2}(0,s;H^{1})}^{2} + \frac{1}{4} \|\sigma\|_{L^{2}(0,s;H^{1})}^{2}  \\
& \quad + \frac{1}{4} \|\sigma_{\infty}\|_{L^{2}(0,T;H^{1})}^{2} \\
& \leq \overline{C} \left ( 2\overline{C} \lambda_{c}^{2} \mathcal{Z}^{2}  + \gamma \varepsilon \lambda_{p} \right ) \|\varphi\|_{L^{2}(0,s;H^{1})}^{2} + \frac{2C_{p}^{2}+1}{4} \|\nabla \sigma\|_{L^{2}(0,s;L^{2})}^{2} \\
& \quad + C \|\sigma_{\infty}\|_{L^{2}(0,T;H^{1})}^{2},
\end{aligned}
\end{equation}
where $\overline{C} := C_{s}^{2} \mathrm{L}_{h} \|\sigma_{1}\|_{L^{\infty}(0,T;L^{2})}$.  Then, upon collecting terms from \eqref{ctsdep:J1J3J4}, \eqref{ctsdep:J2} and \eqref{ctsdep:sigmahvarphi}, we obtain
\begin{equation}\label{ctsdep:truncated:collectionterms}
\begin{aligned}
& \left ( \frac{\gamma \varepsilon}{2} \|\varphi(s)\|_{L^{2}}^{2} +  \mathcal{Z}\frac{\kappa}{4} \|\sigma(s)\|_{L^{2}}^{2} \right ) \\
& \qquad + \int_{0}^{s} \frac{1}{2} \|\mu\|_{L^{2}}^{2} + \left (\frac{1}{2} \mathcal{Z} D_{0} - \frac{2C_{p}^{2} +1}{4} - \frac{2C_{p}^{2}}{4} \left ( 9 \chi^{2} + 7 \right ) \right ) \|\nabla \sigma\|_{L^{2}}^{2} \, dt \\
& \qquad + \int_{0}^{s}  \left ( (\mathcal{Y} - 1) \gamma \varepsilon - 2 \mathcal{Z} D_{0} \eta^{2} - \overline{C} \left ( 2 \overline{C} \lambda_{c}^{2} \mathcal{Z}^{2} + \gamma \varepsilon \lambda_{p} \right )\right ) \|\nabla \varphi\|_{L^{2}}^{2} \, dt \\
& \quad \leq  \frac{\gamma \varepsilon}{2} \|\varphi_{0}\|_{L^{2}}^{2} + \frac{3\kappa}{4}\mathcal{Z}  \|\sigma_{0}\|_{L^{2}}^{2} + \int_{0}^{s}  \mathcal{F} \|\varphi\|_{L^{2}}^{2} \, dt + C\|\mu_{\infty}\|_{L^{2}(0,T;H^{1})}^{2} \\
& \qquad + 2 \kappa \mathcal{Z} \|\sigma_{\infty}\|_{L^{\infty}(0,T;L^{2})}^{2}  + \kappa^{2} \mathcal{Z}^{2} \|\partial_{t} \sigma_{\infty} \|_{L^{2}(0,T;L^{2})}^{2} + C  \| \sigma_{\infty}\|_{L^{2}(0,T;H^{1})}^{2},
\end{aligned}
\end{equation}
where
\begin{align*}
\mathcal{F} & := C \left ( \mathcal{Y} \left ( 1 + \textstyle \sum_{i=1,2} \|\varphi_{i}\|_{H^{3}} \|\varphi_{i}\|_{H^{1}}^{3} \right ) + \left ( 1 + \textstyle \sum_{i=1,2} \|\varphi_{i}\|_{H^{3}}^{2} \|\varphi_{i}\|_{H^{1}}^{6} \right) \right ) \\
& \quad + 3 \mathcal{Y}^{2} + \gamma \varepsilon \lambda_{a} \mathrm{L}_{h} + (\lambda_{p} \gamma \varepsilon h_{\infty})^{2} + \overline{C} \left ( 2 \overline{C} \mathcal{Z}^{2} \lambda_{c}^{2} + \gamma \varepsilon \lambda_{p} \right ).
\end{align*}

We choose
\begin{align*}
\mathcal{Z} > \frac{1}{2D_{0}} \left ( 2C_{p}^{2} \left ( 8 + 9 \chi^{2} \right ) + 1 \right ), \quad \mathcal{Y} > 1 + \frac{1}{\gamma \varepsilon} \left ( 2  \mathcal{Z} D_{0} \eta^{2} + \overline{C} \left (2 \overline{C}  \mathcal{Z}^{2} \lambda_{c}^{2} + \gamma \varepsilon\lambda_{p} \right ) \right )
\end{align*}
so that we obtain from \eqref{ctsdep:truncated:collectionterms} (after adding $ \int_{0}^{s} \mathcal{F} \kappa \|\sigma\|_{L^{2}}^{2} \, dt$ to the right-hand side),
\begin{align*}
& \|\varphi(s)\|_{L^{2}}^{2} + \kappa \|\sigma(s)\|_{L^{2}}^{2} + \int_{0}^{s}  c_{1} \left ( \|\mu\|_{L^{2}}^{2} + \|\nabla \sigma\|_{L^{2}}^{2} + \|\nabla \varphi\|_{L^{2}}^{2} \right ) \, dt \\
& \quad \leq c_{2} \left ( \|\varphi_{0}\|_{L^{2}}^{2} + \kappa \|\sigma_{0}\|_{L^{2}}^{2} \right ) + \int_{0}^{s} c_{3} \left ( 1 + \|\varphi_{i}\|_{H^{3}} + \|\varphi_{i}\|_{H^{3}}^{2} \right ) \left ( \|\varphi\|_{L^{2}}^{2} + \kappa \|\sigma\|_{L^{2}}^{2} \right ) \, dt 
\end{align*}\begin{align*}& \qquad + c_{4} \left ( \|\mu_{\infty}\|_{L^{2}(0,T;H^{1})}^{2} + \|\sigma_{\infty}\|_{L^{2}(0,T;H^{1})}^{2} + \kappa \|\sigma_{\infty}\|_{L^{\infty}(0,T;L^{2})}^{2} \right ) \\
& \qquad + c_{4} \kappa^{2} \|\partial_{t} \sigma_{\infty} \|_{L^{2}(0,T;L^{2})}^{2}.
\end{align*}
for some positive constants $c_{1}$, $c_{2}$, $c_{3}$ and $c_{4}$ independent of $\kappa$ and such that $c_{3}$ depends on $\|\varphi_{i}\|_{L^{\infty}(0,T;H^{1})}$.  Let us define
\begin{align*}
u(s) & := \|\varphi(s)\|_{L^{2}}^{2} + \kappa \|\sigma(s)\|_{L^{2}}^{2}, \quad v(t) := c_{1} \left ( \|\mu\|_{L^{2}}^{2} + \|\nabla \sigma\|_{L^{2}}^{2} + \|\varphi\|_{L^{2}}^{2} \right ), \\
\alpha(t) = \alpha & := c_{2} \left ( \|\varphi_{0}\|_{L^{2}}^{2} + \kappa \|\sigma_{0}\|_{L^{2}}^{2} \right ) + c_{4} \left ( \|\mu_{\infty}\|_{L^{2}(0,T;H^{1})}^{2} + \|\sigma_{\infty}\|_{ L^{2}(0,T;H^{1}) }^{2} \right ) \\
& \quad + c_{4} \left ( \kappa \|\sigma_{\infty}\|_{L^{\infty}(0,T;L^{2})}^{2} + \kappa^{2} \| \partial_{t} \sigma_{\infty} \|_{L^{2}(0,T;L^{2})}^{2} \right ), \\
\beta(t) & := c_{3} \left ( 1 + \|\varphi_{i}(t)\|_{H^{3}} + \|\varphi_{i}(t)\|_{H^{3}}^{2} \right ).
\end{align*}
Thanks to the fact that $\varphi_{i} \in L^{2}(0,T;H^{3})$, $\int_{0}^{T} \beta(t) \, dt$ is finite, and thus
\begin{align*}
\int_{0}^{s} \beta(t) \alpha(t) \exp \left ( \int_{0}^{t} \beta(r) \, dr \right )\, dt
\end{align*}
is also finite for all $s \in (0,T]$.  Then, applying the Gronwall inequality \eqref{Gronwall}, we find that there exists a constant $C$, not depending on $\kappa$, $\varphi$, $\mu$, $\sigma$, $\mu_{\infty}$ and $\sigma_{\infty}$ such that
\begin{align*}
& \|\varphi(s)\|_{L^{2}}^{2}  + \kappa \|\sigma(s)\|_{L^{2}}^{2} + \int_{0}^{s} \left ( \|\mu\|_{L^{2}}^{2} + \|\nabla \sigma\|_{L^{2}}^{2} + \|\nabla \varphi\|_{L^{2}}^{2} \right ) \, dt \\
& \quad \leq C \left ( \|\varphi_{0}\|_{L^{2}}^{2}+ \|\mu_{\infty}\|_{L^{2}(0,T;H^{1})}^{2}  + \|\sigma_{\infty}\|_{L^{2}(0,T;H^{1})}^{2} \right ) \\
& \qquad + C  \kappa \left ( \|\sigma_{0}\|_{L^{2}}^{2} + \|\sigma_{\infty}\|_{L^{\infty}(0,T;L^{2})}^{2} + \kappa \| \partial_{t} \sigma_{\infty} \|_{L^{2}(0,T;L^{2})}^{2} \right )
\end{align*}
for all $s \in (0,T]$.

\vspace*{-5pt}
\section{Quasi-static limit}\label{sec:quasistaticlimit}
\subsection{Existence of weak solutions}
Let $\kappa \in (0,1]$ be fixed and $(\varphi^{\kappa}, \mu^{\kappa}, \sigma^{\kappa})$ be a weak solution to \eqref{Intro:CHN} satisfying the hypothesis of Theorem \ref{thm:truncated:quasi}.  Then, we see that the following energy inequality \eqref{energyineq} is satisfied for $(\varphi^{\kappa}, \mu^{\kappa}, \sigma^{\kappa})$,
\begin{align*}
& \sup_{s \in (0,T]} \left ( \| \Psi(\varphi^{\kappa}(s))\|_{L^{1}} + \|\varphi^{\kappa}(s)\|_{H^{1}}^{2} + \kappa \|\sigma^{\kappa}(s)\|_{L^{2}}^{2} \right ) \\
& \qquad + \|\mu^{\kappa}\|_{L^{2}(0,T;H^{1})}^{2}  + \|\sigma^{\kappa}\|_{L^{2}(0,T;H^{1})}^{2} \\
& \quad \leq C \left ( 1 + \|\partial_{t} \sigma_{\infty}\|_{L^{2}(0,T;L^{2})}^{2} + \|\sigma_{\infty}\|_{L^{\infty}(0,T;L^{2})}^{2} \right ),
\end{align*}
where we have used $\kappa \leq 1$ to obtain that the right-hand side does not depend on $\kappa$.  Together with a similar derivation of \eqref{apriori:time}, this leads to
\begin{align*}
\{\varphi^{\kappa}\}_{\kappa \in (0,1]} & \text{ is bounded in } L^{\infty}(0,T;H^{1}) \cap H^{1}(0,T;H^{-1}), \\
\{\sigma^{\kappa}\}_{\kappa \in (0,1]}, \{\mu^{\kappa}\}_{\kappa \in (0,1]} & \text{ is bounded in } L^{2}(0,T;H^{1}), \\
\{\kappa \partial_{t} \sigma^{\kappa} \}_{\kappa \in (0,1]} & \text{ is bounded in } L^{2}(0,T;H^{-1}),
\end{align*}
and by the standard compactness results, there exist functions $\varphi_{*}$, $\mu_{*}$ and $\sigma_{*}$ such that, for a non-relabeled subsequence,
\begin{alignat*}{3}
\varphi^{\kappa} & \rightarrow \varphi_{*} && \quad \text{ weakly-}* && \quad \text{ in } L^{\infty}(0,T;H^{1}) \cap H^{1}(0,T;H^{-1}), \\
\varphi^{\kappa} & \rightarrow \varphi_{*} && \quad \text{ strongly } && \quad \text{ in } C^{0}([0,T];H^{r}) \cap L^{2}(0,T;L^{m}) \text{ and a.e. in } \Omega \times (0,T), 
\end{alignat*}\begin{alignat*}{3}\sigma^{\kappa} & \rightarrow \sigma_{*} && \quad \text{ weakly } && \quad \text{ in } L^{2}(0,T;H^{1}), \\
\partial_{t}(\kappa \sigma^{\kappa}) & \rightarrow 0 && \quad \text{ weakly } && \quad \text{ in } L^{2}(0,T;H^{-1}), \\
\mu^{\kappa} & \rightarrow \mu_{*} && \quad \text{ weakly } && \quad \text{ in } L^{2}(0,T;H^{1}),
\end{alignat*}
where $0 \leq r < 1$ and $m \in [1,6)$.  Then, passing to the limit $\kappa \to 0$ in \eqref{truncated:weak} shows that $(\varphi_{*}, \mu_{*}, \sigma_{*})$ is a weak solution triplet to \eqref{Intro:quasistatic} in the sense of Definition \ref{defn:quasi:weak}.

\subsection{Further regularity}
The $L^{2}(0,T;H^{3})$-regularity for $\varphi_{*}$ follows along the same lines as in Section \ref{sec:furtherreg} and so we will omit the details.  For the $L^{\infty}(0,T;H^{1})$-regularity for $\sigma_{*}$, we turn to \eqref{quasistatic:weak:sigma} and test with $\xi = \sigma_{*} - \sigma_{\infty}$, resulting in
\begin{align*}
&  \int_{\Omega} D(\varphi_{*}) |\nabla \sigma_{*}|^{2} + \lambda_{c} h(\varphi_{*}) |\sigma_{*}|^{2} \, dx \\
& \quad = \int_{\Omega} D(\varphi_{*}) \left ( \eta \nabla \varphi_{*} \cdot \nabla (\sigma_{*} - \sigma_{\infty}) + \nabla \sigma_{*} \cdot \nabla \sigma_{\infty} \right ) + \lambda_{c} h(\varphi_{*}) \sigma_{*} \sigma_{\infty} \, dx.
\end{align*}
By Young's inequality, and the Poincar\'{e} inequality, we find that
\begin{align*}
\|\nabla \sigma_{*}\|_{L^{2}}^{2} \leq C(C_{p}, D_{0}, D_{1}, \eta, \lambda_{c}, h_{\infty}) \left ( \|\sigma_{\infty}\|_{H^{1}}^{2} + \|\nabla \varphi_{*}\|_{L^{2}}^{2} \right ),
\end{align*}
where the right-hand side is bounded in $L^{\infty}(0,T)$.  Then, by the Poincar\'{e} inequality, we easily have
\begin{align*}
\sigma_{*} \in L^{\infty}(0,T;H^{1}).
\end{align*}

\subsection{Continuous dependence}
The continuous dependence result can be easily obtain by setting $\kappa = 0$ in Section \ref{sec:ctsdep}.  Given two weak solution triplets satisfying the hypothesis of Theorem \ref{thm:truncated:quasi:ctsdep}, let $\varphi$, $\mu$, and $\sigma$ denote their differences, respectively, which satisfy
\begin{align*}
0 & = \langle \partial_{t}\varphi, \zeta \rangle + \int_{\Omega} \nabla \mu \cdot \nabla \zeta + \left ( \lambda_{a} - \lambda_{p} \sigma_{1} \right ) ( h(\varphi_{1}) - h(\varphi_{2})) \zeta - \lambda_{p} \sigma h(\varphi_{2}) \zeta \, dx, \\
0 & = \int_{\Omega} \mu \lambda - \frac{\gamma}{\varepsilon} (\Psi'(\varphi_{1}) - \Psi'(\varphi_{2})) \lambda - \gamma \varepsilon \nabla \varphi \cdot \nabla \lambda + \chi \sigma \lambda \, dx, \\
0 & = \int_{\Omega} D_{0}(\nabla \sigma - \eta \nabla \varphi) \cdot \nabla \xi + \lambda_{c}( \sigma_{1}( h(\varphi_{1}) - h(\varphi_{2})) + \sigma h(\varphi_{2})) \xi \, dx,
\end{align*}
for all $\zeta, \lambda, \xi \in H^{1}_{0}(\Omega)$.  As in Section \ref{sec:ctsdep}, let $\mu_{\infty} := \mu_{\infty,1} - \mu_{\infty,2}$ and $\sigma_{\infty} := \sigma_{\infty,1} - \sigma_{\infty,2}$ denote the difference of boundary data, and we substitute $\zeta = \gamma \varepsilon \varphi$, $\xi = \mathcal{Z} (\sigma - \sigma_{\infty})$, $\lambda = \mu - \mu_{\infty}$ and $\lambda = \mathcal{Y} \varphi$ for positive constants $\mathcal{Y}$, $\mathcal{Z}$ yet to be determined.  This is equivalent to setting $\kappa = 0$ in \eqref{CHN:uniqueness:apriori:est}.  Then, $J_{1}$ vanishes and the estimations of $J_{2}$, $J_{3}$, $J_{4}$ and $J_{5}$ do not depend on $\kappa$.  One would have to adjust the value for $\mathcal{Z}$ accordingly to account for the absence of the first term on the right-hand side of \eqref{ctsdep:J1}, but the same arguments will lead to the assertion of Theorem \ref{thm:truncated:quasi:ctsdep}.

\section{Singular potentials}\label{sec:singularPot}

\subsection{Maximal monotone operators and the Yosida approximation}\label{sec:Yosida}
In this section, we will briefly review the basic concepts regarding the Yosida approximation of maximal monotone operators.  For more details related to maximal monotone operators, subdifferentials and the Yosida approximation, we refer the reader to \cite[Chapter III]{Brezis}, \cite[Chapter III]{Pascali}, \cite[p. 161]{Showalter} and \cite[Chapter 32]{Zeidler}.

Let $X$ be a Banach space with dual space $X^{*}$.  A possibly multivalued mapping $T : X \to 2^{X^{*}}$ with effective domain $D(T) := \{ x \in X : Tx \neq \emptyset \}$, range $R(T) := \{y \in Tx : x \in D(T) \}$ and graph $G(T) := \{ [x,y] \in X \times X^{*} : x \in D(T), y \in Tx \}$ is said to be monotone if
\begin{align*}
\langle f-g, x-y \rangle \geq 0 \quad \forall x,y \in D(T), \; f \in Tx, \; g \in Ty.
\end{align*}
The mapping $T$ is maximal monotone if and only if
\begin{align*}
\langle f-g, x-y \rangle \geq 0 \quad \forall [y,g] \in G(T) \Rightarrow x \in D(T) \text{ and } f \in Tx.
\end{align*}
In the case where $X$ is a Hilbert space $H$, we identify $H$ with its dual $H^{*}$.  For a maximal monotone multivalued mapping $T: H \to 2^{H}$ with $[0,0] \in G(T)$, i.e., $0 \in T(0)$, and for any $\lambda > 0$, one can show that $I + \lambda T$ is surjective, where $I$ is the identity operator, i.e., $R(I + \lambda T) = H$.  Hence, for every fixed $u \in H$, there exists $[x,f] \in G(T)$ such that
\begin{align*}
u = x + \lambda f \text{ with } f \in Tx.
\end{align*}
The resolvent $\mathcal{J}_{\lambda} : H \to D(T)$ and the Yosida approximant $T_{\lambda} : H \to H$ of $T$ are defined as
\begin{align*}
& \mathcal{J}_{\lambda} := (I + \lambda T)^{-1}, \quad T_{\lambda} := \frac{1}{\lambda} (I - \mathcal{J}_{\lambda}), \text{ so that } x = \mathcal{J}_{\lambda} u, \; f = T_{\lambda} u \\
& \text{ satisfy } u = \mathcal{J}_{\lambda}u + \lambda T_{\lambda} u \quad \forall \lambda > 0 \text{ and } T_{\lambda} u = f \in Tx = T \left ( \mathcal{J}_{\lambda}u \right ).
\end{align*}
It is well-known that $\mathcal{J}_{\lambda}$ and $T_{\lambda}$ are both Lipschitz operators with Lipschitz constant $1$ and $\frac{1}{\lambda}$, respectively.  Furthermore, denoting $T^{0}x$ as the element such that $\|T^{0}x\| := \inf \{ \|f\|_{H} : f \in Tx \}$, then it holds that
\begin{align}\label{Yosidaweakconv}
\|T_{\lambda}x\|_{H} \leq \|T^{0}x\| \quad \forall \lambda > 0 \text{ and } T_{\lambda} x \rightharpoonup T^{0}x \text{ in } H \text{ as } \lambda \to 0 \quad \forall x \in D(T).
\end{align}
For a convex, proper, lower semicontinuous function $\hat{\beta}: H \to \overline{\mathbb{R}}$, it is also well-known that the subdifferential $T = \partial \hat{\beta} : H \to 2^{H}$, defined for $x \in D(\partial \hat{\beta}) = \{ y \in H : \partial \hat{\beta}(y) \neq \emptyset \}$ as
\begin{align*}
\partial \hat{\beta}(x) := \{ f \in H : \hat{\beta}(y) - \hat{\beta}(x) \geq \langle f, y-x \rangle \; \forall y \in H \}
\end{align*}
is a maximal monotone mapping.  For any $\lambda > 0$, the Yosida approximation of $\hat{\beta}$, defined as
\begin{align*}
\hat{\beta}_{\lambda}(x) = \inf_{y \in H} \left \{ \frac{1}{2\lambda} \|x - y\|_{H}^{2} + \hat{\beta}(y)\right \},
\end{align*}
is convex, Fr\'{e}chet differentiable on $H$ with $\mathrm{d} \hat{\beta}_{\lambda} = \partial \hat{\beta}_{\lambda} = T_{\lambda} = (\partial \hat{\beta})_{\lambda}$ for all $\lambda > 0$, i.e., the Fr\'{e}chet derivative $\mathrm{d} \hat{\beta}_{\lambda}$ of $\hat{\beta}_{\lambda}$ coincides with its subdifferential $\partial \hat{\beta}_{\lambda}$ and is equal to the Yosida approximant of the subdifferential $\partial \hat{\beta}$.  Moreover, it holds that
\begin{align}\label{Yosida:monotoneconvUP}
\hat{\beta}_{\lambda}(x) \leq \hat{\beta}(x) \text{ for all } x \in H \text{ and } \hat{\beta}_{\lambda}(x) \nearrow \hat{\beta}(x) \text{ as } \lambda \to 0.
\end{align}

\subsection{Existence}
From the above discussion, we take the Hilbert space $H = \mathbb{R}$ and introduce the Yosida approximation of $\beta = \partial \hat{\beta}$ as follows.  For $n \in (0,1]$, let
\begin{align}
\label{Yosida}
\hat{\beta}_{n}(y) := \min_{s \in \mathbb{R}} \left ( \frac{1}{2n} (s-y)^{2} + \hat{\beta}(s) \right ) \geq 0, \quad \beta_{n}(y) := \hat{\beta}_{n}'(y),
\end{align}
and
\begin{align*}
\Psi_{n}(y) := \hat{\beta}_{n}(y) + \Lambda(y).
\end{align*}
Then, for each $n \in (0,1]$, $\Psi_{n} \in C^{1,1}(\mathbb{R})$ is non-negative and has at most quadratic growth.  The conditions of Theorems \ref{thm:truncated:exist} and \ref{thm:truncated:regularity} are fulfilled (with $p = 1)$, and so, for every $n \in (0,1]$, there exists a weak solution $(\varphi_{n}, \mu_{n}, \sigma_{n})$ to \eqref{Intro:CHN} with $\Psi'_{n}$ in \eqref{Intro:mu} and satisfies $\varphi_{n}(0) = \varphi_{0}$, $\sigma_{n}(0) = \sigma_{0}$ with the following regularity
\begin{align*}
\varphi_{n} & \in (-1 + L^{\infty}(0,T;H^{1}_{0})) \cap L^{2}(0,T;H^{3}) \cap H^{1}(0,T;H^{-1}), \\
\mu_{n} & \in \mu_{\infty} + L^{2}(0,T;H^{1}_{0}), \\
\sigma_{n} & \in (\sigma_{\infty} + L^{2}(0,T;H^{1}_{0})) \cap L^{\infty}(0,T;L^{2}) \cap H^{1}(0,T;H^{-1}).
\end{align*}
Moreover, the weak solution $(\varphi_{n}, \mu_{n}, \sigma_{n})$ satisfies an analogue of \eqref{energyineq}:
\begin{equation}\label{energyineq:Yosida}
\begin{aligned}
& \sup_{s \in (0,T]} \left ( \|\Psi_{n}(\varphi_{n}(s))\|_{L^{1}} + \|\varphi_{n}(s)\|_{H^{1}}^{2} + \kappa \|\sigma_{n}(s)\|_{L^{2}}^{2} \right ) \\
& \qquad  + \|\mu_{n}\|_{L^{2}(0,T;H^{1})}^{2} + \|\sigma_{n}\|_{L^{2}(0,T;H^{1})}^{2} \\
& \quad \leq  C \left ( 1 + \kappa \| \sigma_{0} - \sigma_{\infty}(0)\|_{L^{2}}^{2} + \kappa^{2} \|\partial_{t} \sigma_{\infty}\|_{L^{2}(0,T;L^{2})}^{2} + \kappa \|\sigma_{\infty}\|_{L^{\infty}(0,T;L^{2})}^{2} \right ),
\end{aligned}
\end{equation}
for some positive constant $C$ independent of $n$ and $\kappa$.  The fact that $C$ is independent of $n$ follows from the fact that that the derivation of \eqref{energyineq} does not use the polynomial growth of the potential.  Moreover, by (S3) and \eqref{Yosida:monotoneconvUP}, the quantity $\|\Psi_{n}(\varphi_{n}(0))\|_{L^{1}} = \|\Psi_{n}(\varphi_{0})\|_{L^{1}}$ is bounded uniformly in $n$.  This allows us to deduce that
\begin{align*}
\{\varphi_{n}\}_{n \in (0,1]} & \text{ is bounded in } L^{\infty}(0,T;H^{1}), \\
\{\mu_{n}\}_{n \in (0,1]} & \text{ is bounded in } L^{2}(0,T;H^{1}), \\
\{\sigma_{n}\}_{n \in (0,1]} & \text{ is bounded in } L^{\infty}(0,T;L^{2}) \cap L^{2}(0,T;H^{1}).
\end{align*}
By the Lipschitz continuity of $\Lambda'$, we also obtain that $\{\Lambda'(\varphi_{n})\}_{n \in (0,1]}$ is bounded in $L^{\infty}(0,T;H^{1})$.  Next, due to the regularity of $(\varphi_{n}, \mu_{n}, \sigma_{n})$, we can write the equation for $\mu_{n}$ as the following equality
\begin{align}\label{SingularPotential:munequ}
 \gamma \varepsilon \Delta \varphi_{n}  = - \mu_{n} - \chi \sigma_{n} + \frac{\gamma}{\varepsilon} \Lambda_{n}'(\varphi_{n}) + \frac{\gamma}{\varepsilon} \beta_{n}(\varphi_{n}).
\end{align}
As $\beta_{n}$ is Lipschitz and therefore belong to $W^{1,\infty}(\mathbb{R})$, it is differentiable with derivative $\beta_{n}' \in L^{\infty}(\mathbb{R})$ a.e. in $\mathbb{R}$.  Then, by the chain rule for Sobolev functions \cite[Thm 2.1.11]{Ziemer}, it holds that
\begin{align*}
\int_{\Omega} - \Delta \varphi_{n} (\beta_{n}(\varphi_{n}) - \beta_{n}(-1)) \, dx = \int_{\Omega} \beta_{n}'(\varphi_{n}) |\nabla \varphi_{n}|^{2} \, dx,
\end{align*}
where we used that $\beta_{n}(\varphi_{n}) - \beta_{n}(-1) = 0$ on $\Gamma$.  Here we point out that this is where the assumption $-1 \in D(\beta)$ comes in, as we can only test with $H^{1}_{0}$ functions for \eqref{SingularPotential:munequ}.  Thus, after multiplying \eqref{SingularPotential:munequ} by $\beta_{n}(\varphi_{n}) - \beta_{n}(-1) \in H^{1}_{0}$ and performing integration by parts, we obtain
\begin{align*}
& \int_{\Omega} \frac{\gamma}{\varepsilon} |\beta_{n}(\varphi_{n})|^{2} + \beta_{n}'(\varphi_{n}) |\nabla \varphi_{n}|^{2} \, dx \\
& \quad = \int_{\Omega} \left ( -\frac{\gamma}{\varepsilon} \Lambda'(\varphi_{n}) + (\mu_{n} + \chi \sigma_{n}) \right ) (\beta_{n}(\varphi_{n}) - \beta_{n}(-1)) + \frac{\gamma}{\varepsilon} \beta_{n}(\varphi_{n}) \beta_{n}(-1) \, dx.
\end{align*}

Using that $\hat{\beta}_{n}$ is convex and thus $\beta_{n}' \geq 0$ a.e. in $\mathbb{R}$, the second term on the left-hand side can be neglected.  Furthermore by Assumption \ref{assump:singularpotential} and \eqref{Yosidaweakconv}, we have $|\beta_{n}(-1)| \leq |\beta^{0}(-1)| < \infty$ for all $n \in (0,1]$.  Thus, writing $\Lambda'
(\varphi_{n}) = \Lambda'(\varphi_{n}) - \Lambda'(0) + \Lambda'(0)$, and using the Lipschitz continuity of $\Lambda'$, H\"{o}lder's inequality and Young's inequality, we have
\begin{align*}
\|\beta_{n}(\varphi_{n})\|_{L^{2}(0,T;L^{2})}^{2} \leq C \left ( 1 +  \|\varphi_{n}\|_{L^{2}(0,T;L^{2})}^{2} + \|\mu_{n} + \chi \sigma_{n}\|_{L^{2}(0,T;L^{2})}^{2}  \right ),
\end{align*}
for some positive constant $C$ depending only on $\gamma$, $\varepsilon$, $\|\Lambda''\|_{L^{\infty}(\mathbb{R})}$, $|\Lambda'(0)|$, and $|\beta^{0}(-1)|$.  Returning to \eqref{SingularPotential:munequ}, as the right-hand side is now bounded in $L^{2}(0,T;L^{2})$ uniformly in $n$ and $\kappa$, by elliptic regularity, we have that $\{\varphi_{n}\}_{n \in (0,1]}$ is also bounded in $L^{2}(0,T;H^{2})$.

Analogous to \eqref{apriori:time}, $\{ \partial_{t}\varphi_{n}\}_{n \in (0,1]}$ and $\{ \kappa \partial_{t} \sigma_{n}\}_{n \in (0,1]}$ are bounded in $L^{2}(0,T;H^{-1})$.  Hence, in addition to the convergence stated in Section \ref{sec:passinglimit}, where we denote the limit functions of $(\varphi_{n}, \mu_{n}, \sigma_{n})$ as $(\varphi, \mu, \sigma)$ (changing the index from $k$ to $n$), we also have
\begin{alignat*}{3}
\varphi_{n} & \rightarrow \varphi && \quad \text{ weakly } && \quad \text{ in } L^{2}(0,T;H^{2}), \\
\beta_{n}(\varphi_{n}) & \rightarrow \psi && \quad \text{ weakly } && \quad \text{ in } L^{2}(0,T;L^{2})
\end{alignat*}
to some function $\psi \in L^{2}(0,T;L^{2})$.

To finish the proof, it suffices to pass to the limit $n \to 0$ in the weak formulation of $(\varphi_{n}, \mu_{n}, \sigma_{n})$ to show that $(\varphi, \mu, \sigma, \psi)$ satisfies \eqref{singular:weak}.  We will omit the details.  It remains to show that $\varphi \in D(\beta)$ and $\psi \in \beta(\varphi)$ a.e. in $\Omega \times (0,T)$, which will follow from the maximal monotonicity of $\beta$ once we showed
\begin{align*}
\langle \psi - g , \varphi - y \rangle \geq 0 \quad \forall [y,g] \in G(\beta).
\end{align*}
We argue as in \cite[Lem 1.3 (e), p. 127]{BCP70} (cf. \cite[Prop 2.2 (iv), p. 38]{Barbu}).  By the boundedness of $\{\beta_{n}(\varphi_{n})\}_{n \in (0,1]}$ in $L^{2}(0,T;L^{2})$ and the identity
\begin{align*}
\varphi_{n} = \mathcal{J}_{n} \varphi_{n} + n \beta_{n} (\varphi_{n}),
\end{align*}
where $\mathcal{J}_{n} \varphi_{n}$ is the resolvent of $\varphi_{n}$, we find that $\varphi_{n} - \mathcal{J}_{n}\varphi_{n} \to 0$ strongly in $L^{2}(0,T;L^{2})$, and this implies that $\mathcal{J}_{n}\varphi_{n} \to \varphi$ strongly in $L^{2}(0,T;L^{2})$.  Take an arbitrary $[y,g] \in G(\beta)$, i.e., $y \in D(\beta)$ and $g \in \beta(y)$.  By the monotonicity of $\beta$ and as $\beta_{n}(\varphi_{n}) \in \beta (\mathcal{J}_{n} \varphi_{n})$, we have (here $\langle \cdot, \cdot \rangle$ denotes the scalar product on $L^{2}(0,T;L^{2})$)
\begin{align*}
0 \leq \langle \beta_{n}(\varphi_{n}) - g, \mathcal{J}_{n}\varphi_{n} - y \rangle \to \langle \psi - g, \varphi - y \rangle \text{ as } n \to \infty.
\end{align*}
As $[y,g] \in G(\beta)$ is arbitrary, by the maximal monotonicity of $\beta$, we have that $[\varphi, \psi] \in G(\beta)$, and so $\varphi \in D(\beta)$ and $\psi \in \beta(\varphi)$ a.e. in $\Omega \times (0,T)$.

\subsection{Energy inequality}
We argue as in Section \ref{sec:energyineq} and pass to the limit $n \to 0$ in \eqref{energyineq:Yosida}.  Due to the non-negativity and continuity of $\Lambda(\cdot)$, and the a.e. convergence $\varphi_{n} \to \varphi$ in $\Omega \times (0,T)$, by Fatou's lemma we have for a.e. $s \in (0,T]$,
\begin{align*}
\int_{\Omega} \Lambda(\varphi(s)) \, dx & \leq \liminf_{n \to 0} \int_{\Omega} \Lambda(\varphi_{n}(s)) \, dx.
\end{align*}
Similarly, for every $n \in (0,1]$, $\hat{\beta}_{n}$ is non-negative and continuous.  By \eqref{Yosida:monotoneconvUP} and a.e. convergence $\varphi_{n} \to \varphi$ in $\Omega \times (0,T)$, we have $\hat{\beta}_{n}(\varphi_{n}) \to \hat{\beta}(\varphi)$ a.e. in $\Omega \times (0,T)$.  By Fatou's lemma, we obtain for a.e. $s \in (0,T]$,
\begin{align*}
\int_{\Omega} \hat{\beta}(\varphi(s)) \, dx & \leq \liminf_{n \to 0} \int_{\Omega} \hat{\beta}_{n}(\varphi_{n}(s)) \, dx.
\end{align*}
Using weak/weak-$*$ lower semicontinuity of the Sobolev norms on the other terms, and passing to the limit $n \to 0$ in \eqref{energyineq:Yosida} leads to \eqref{energyineq:singular}.

\subsection{Partial continuous dependence}
Consider two weak solution quadruples $(\varphi_{i}, \mu_{i}, \sigma_{i}, \psi_{i})_{i = 1,2}$ to \eqref{singular:CHN} satisfying the hypothesis of Theorem \ref{thm:singular:ctsdep}.  Denoting the differences of $(\varphi_{i}, \mu_{i}, \sigma_{i}, \psi_{i}, \sigma_{\infty,i})_{i=1,2}$ as $\varphi$, $\mu$, $\sigma$, $\psi$ and $\sigma_{\infty}$, respectively, we then have
\begin{align*}
\varphi & \in L^{\infty}(0,T;H^{1}_{0}) \cap L^{2}(0,T;H^{2}) \cap H^{1}(0,T;H^{-1}), \\
\sigma & \in \left ( \sigma_{\infty} +  L^{2}(0,T;H^{1}_{0}) \right ) \cap L^{\infty}(0,T;L^{2}) \cap H^{1}(0,T;H^{-1}), \\
\mu & \in L^{2}(0,T;H^{1}_{0} ), \quad \psi \in L^{2}(0,T;L^{2}),
\end{align*}
and
\begin{align}
\label{singular:truncated:uniq:varphi} 0 & = \langle \partial_{t}\varphi, \zeta \rangle + \int_{\Omega} \nabla \mu \cdot \nabla \zeta + (h(\varphi_{1}) - h(\varphi_{2})) (\lambda_{a} - \lambda_{p} \sigma_{1}) \zeta - \lambda_{p} \sigma h(\varphi_{2}) \zeta \,dx, \\
\label{singular:truncated:uniq:mu} 0 & = \int_{\Omega} \mu \lambda - \frac{\gamma}{\varepsilon} (\psi + \Lambda'(\varphi_{1}) - \Lambda'(\varphi_{2})) \lambda - \gamma \varepsilon \nabla \varphi \cdot \nabla \lambda + \chi \sigma \lambda \,dx, \\
\notag 0 & = \langle \kappa \partial_{t} \sigma, \xi \rangle + \int_{\Omega} D_{0}(\nabla \sigma - \eta \nabla \varphi) \cdot \nabla \xi + \lambda_{c} (\sigma_{1} (h(\varphi_{1}) - h(\varphi_{2})) + \sigma h(\varphi_{2})) \xi \, dx,
\end{align}
for all $\zeta$, $\lambda$, $\xi \in H^{1}_{0}$ and a.e. $t \in (0,T)$.

Let $\mathcal{Y}$ denote a positive constant yet to be determined.  Testing with $\zeta = \mathcal{Y} N(\varphi) \in H^{1}_{0}$, where $N$ is the inverse Dirichlet-Laplacian operator defined in \eqref{InverseLaplacian}, and $\lambda = - \mathcal{Y}\varphi \in H^{1}_{0}$ leads to
\begin{align*}
0 & = \frac{\mathcal{Y}}{2} \frac{d}{dt} \|\varphi\|_{*}^{2} + \mathcal{Y} \int_{\Omega} \nabla \mu \cdot \nabla N(\varphi) \, dx \\
& \quad + \mathcal{Y} \int_{\Omega}  ((\lambda_{a} - \lambda_{p} \sigma_{1})(h(\varphi_{1})-h(\varphi_{2})) - \lambda_{p} \sigma h(\varphi_{2})) N(\varphi) \, dx, \\
0 & = \mathcal{Y} \int_{\Omega} \frac{\gamma}{\varepsilon} (\psi \varphi + (\Lambda'(\varphi_{1}) - \Lambda'(\varphi_{2})) \varphi + \gamma \varepsilon |\nabla \varphi|^{2} - (\mu + \chi \sigma) \varphi \, dx,
\end{align*}
where we used \eqref{TimederivativeStar}.  Upon adding the above equations leads to
\begin{align*}
\frac{\mathcal{Y}}{2} \frac{d}{dt} \|\varphi\|_{*}^{2} &  + \mathcal{Y} \int_{\Omega} \nabla \mu \cdot \nabla N(\varphi) - \mu \varphi + \frac{\gamma}{\varepsilon} \psi \varphi \, dx \\
& + \mathcal{Y} \int_{\Omega} ((\lambda_{a} - \lambda_{p} \sigma_{1})(h(\varphi_{1})-h(\varphi_{2})) - \lambda_{p} \sigma h(\varphi_{2})) N(\varphi) \, dx \\
& + \mathcal{Y} \int_{\Omega} \frac{\gamma}{\varepsilon} ( \Lambda'(\varphi_{1}) - \Lambda'(\varphi_{2})) \varphi + \gamma \varepsilon |\nabla \varphi|^{2} - \chi \sigma \varphi \, dx = 0.
\end{align*}
Using \eqref{InverseLaplacian} with $u = \mu \in H^{1}_{0}$, we observe that
\begin{align*}
\int_{\Omega} \nabla \mu \cdot \nabla N(\varphi) \,dx - \int_{\Omega} \mu \varphi \, dx = 0,
\end{align*}
and by the monotonicity of $\beta$, for $\psi_{i} \in \beta(\varphi_{i})$, $i = 1,2$, we have
\begin{align*}
\int_{\Omega} \psi \varphi \, dx = \int_{\Omega} (\psi_{1} - \psi_{2})(\varphi_{1} - \varphi_{2}) \, dx \geq 0.
\end{align*}
Thus, by H\"{o}lder's inequality and the Lipschitz continuity of $\Lambda'$, it holds that
\begin{equation}
\label{ctsdep:singular:1}
\begin{aligned}
 & \frac{\mathcal{Y}}{2} \frac{d}{dt} \|\varphi\|_{*}^{2} + \mathcal{Y} \gamma \varepsilon \|\nabla \varphi\|_{L^{2}}^{2} \\
& \quad \leq \mathcal{Y} \left ( \lambda_{a} \mathrm{L}_{h} \|\varphi\|_{L^{2}} \|N(\varphi)\|_{L^{2}} + \lambda_{p} \mathrm{L}_{h} \|\sigma_{1}\|_{L^{2}} \|\varphi\|_{L^{4}} \|N(\varphi)\|_{L^{4}} \right ) \\
& \qquad + \mathcal{Y} \left ( \lambda_{p} h_{\infty} \|\sigma\|_{L^{2}} \|N(\varphi)\|_{L^{2}} + \frac{\gamma}{\varepsilon} \|\Lambda''\|_{L^{\infty}(\mathbb{R})} \|\varphi\|_{L^{2}}^{2} + \chi \|\sigma\|_{L^{2}} \|\varphi\|_{L^{2}} \right ) .
\end{aligned}
\end{equation}
Adding \eqref{uniqueness:sigmaequ:testedsigma} to \eqref{ctsdep:singular:1} we obtain (neglecting the non-negative term $\mathcal{Z} \lambda_{c} \int_{\Omega} h(\varphi_{2}) |\sigma|^{2}\break \, dx$) after integrating in time from $0$ to $s \in (0,T]$,
\begin{equation}\label{ctsdep:singular:combined}
\begin{aligned}
& \frac{\mathcal{Y}}{2} \|\varphi(s)\|_{*}^{2} + \mathcal{Z} \frac{\kappa}{2} \|\sigma(s)\|_{L^{2}}^{2}  \\
& \qquad + \mathcal{Y} \gamma \varepsilon \|\nabla \varphi\|_{L^{2}(0,s;L^{2})}^{2} + \mathcal{Z} D_{0} \|\nabla \sigma\|_{L^{2}(0,s;L^{2})}^{2} \\
& \quad \leq J_{1} + K_{1} + K_{2} + K_{3} + K_{4} + \frac{\mathcal{Y}}{2} \|\varphi_{0}\|_{*}^{2} + \mathcal{Z} \frac{\kappa}{2} \|\sigma_{0}\|_{L^{2}}^{2},
\end{aligned}
\end{equation}
where $J_{1}$, as defined in Section \ref{sec:ctsdep}, can be handled as in \eqref{ctsdep:J1}, and
\begin{align*}
K_{1} & = \mathcal{Y} \int_{0}^{s} \left ( \lambda_{a} \mathrm{L}_{h} \|\varphi\|_{L^{2}} \|N(\varphi)\|_{L^{2}} + \lambda_{p} \mathrm{L}_{h} \|\sigma_{1}\|_{L^{2}} \|\varphi\|_{L^{4}} \|N(\varphi)\|_{L^{4}} \right ) \, dt, \\
K_{2} & = \mathcal{Y} \int_{0}^{s} \left ( \lambda_{p} h_{\infty} \|\sigma\|_{L^{2}} \|N(\varphi)\|_{L^{2}} + \frac{\gamma}{\varepsilon} \|\Lambda''\|_{L^{\infty}(\mathbb{R})} \|\varphi\|_{L^{2}}^{2} + \chi \|\sigma\|_{L^{2}} \|\varphi\|_{L^{2}} \right ) \, dt, \\
K_{3} & = \mathcal{Z} \int_{0}^{s} \int_{\Omega} D_{0} \eta \nabla \varphi \cdot \nabla (\sigma - \sigma_{\infty}) + D_{0} \nabla \sigma \cdot \nabla \sigma_{\infty} \, dx \, dt, \\
K_{4} & =- \mathcal{Z} \int_{0}^{s} \int_{\Omega} \lambda_{c} \sigma_{1} (h(\varphi_{1}) - h(\varphi_{2}))(\sigma - \sigma_{\infty}) - \lambda_{c} h(\varphi_{2}) \sigma \sigma_{\infty} \, dx \, dt.
\end{align*}
By H\"{o}lder's inequality and Young's inequality, and also from the estimation of $J_{5}$ in \eqref{ctsdep:sigmahvarphi}, we see that
\begin{subequations}
\begin{align}
\label{singular:ctsdep:K3} |K_{3}| & \leq  \mathcal{Z} \frac{2D_{0}}{4} \|\nabla \sigma\|_{L^{2}(0,s;L^{2})}^{2} + \mathcal{Z} \frac{5D_{0}}{4} \|\nabla \sigma_{\infty}\|_{L^{2}(0,T;L^{2})}^{2} \\
\notag & \quad + 2 \mathcal{Z} D_{0} \eta^{2} \|\nabla \varphi\|_{L^{2}(0,s;L^{2})}^{2}, \\
\label{singular:ctsdep:K4} |K_{4}| & \leq 2 \overline{C}^{2} \lambda_{c}^{2} \mathcal{Z}^{2} \|\varphi\|_{L^{2}(0,s;H^{1})}^{2} + \frac{2}{4} \|\sigma\|_{L^{2}(0,s;H^{1})}^{2} \\
\notag & \quad + \left ( \frac{1}{4} + (\lambda_{c} h_{\infty} \mathcal{Z})^{2} \right ) \|\sigma_{\infty}\|_{L^{2}(0,T;H^{1})}^{2} \\
\notag & \leq 2 \overline{C}^{2} \lambda_{c}^{2} \mathcal{Z}^{2} (C_{p}+1)^{2} \|\nabla \varphi\|_{L^{2}(0,s;L^{2})}^{2} + \frac{2C_{p}^{2} + 1}{2} \|\nabla \sigma\|_{L^{2}(0,s;L^{2})}^{2} \\
\notag & \quad + C \|\sigma_{\infty}\|_{L^{2}(0,T;H^{1})}^{2},
\end{align}
\end{subequations}
where $\overline{C} = C_{s}^{2} \mathrm{L}_{h} \|\sigma_{1}\|_{L^{\infty}(0,T;L^{2})}$ as before, and we have applied the Poincar\'{e} inequality to $\varphi \in H^{1}_{0}$.  Meanwhile, using the Poincar\'{e} inequality applied to $N(\varphi) \in H^{1}_{0}$, we see that
\begin{align*}
\|N(\varphi)\|_{H^{1}} \leq (C_{p}+1) \|\nabla N(\varphi)\|_{L^{2}} = (C_{p}+1) \|\varphi\|_{*},
\end{align*}
and by the Sobolev embedding $H^{1} \subset L^{4}$, it holds that
\begin{align*}
\|\varphi\|_{L^{4}} \|N(\varphi)\|_{L^{4}} & \leq C_{s}^{2} \|\varphi\|_{H^{1}} \|N(\varphi)\|_{H^{1}} \leq C_{s}^{2} (C_{p}+1)^{2} \|\varphi\|_{*} \|\nabla \varphi\|_{L^{2}}.
\end{align*}
Then, by H\"{o}lder's inequality and Young's inequality it holds that
\begin{align*}
|K_{1}| & \leq (\mathcal{Y} \lambda_{a} \mathrm{L}_{h} C_{p} )^{2} \|N(\varphi)\|_{L^{2}(0,s;L^{2})}^{2} + \frac{1}{4} \|\nabla \varphi\|_{L^{2}(0,s;L^{2})}^{2} \\
 & \quad + (\mathcal{Y} \lambda_{p} \overline{C} (C_{p}+1)^{2})^{2} \|\varphi\|_{L^{2}(0,s;*)}^{2} + \frac{1}{4} \|\nabla \varphi\|_{L^{2}(0,s;L^{2})}^{2}, \\
|K_{2}| & \leq (\mathcal{Y} \lambda_{p} h_{\infty})^{2} \|N(\varphi)\|_{L^{2}(0,s;L^{2})}^{2}  + \frac{1}{4} \|\sigma\|_{L^{2}(0,s;L^{2})}^{2} \\
& \quad + \left ( \mathcal{Y} \frac{\gamma}{\varepsilon} \|\Lambda''\|_{L^{\infty}(\mathbb{R})} + \frac{\mathcal{Y}^{2}}{2} \right ) \|\varphi\|_{L^{2}(0,s;L^{2})}^{2} + \frac{\chi^{2}}{2} \|\sigma\|_{L^{2}(0,s;L^{2})}^{2},
\end{align*}
where we use the notation
\begin{align*}
\|f\|_{L^{2}(0,s;*)}^{2} = \int_{0}^{s} \|f\|_{*}^{2} \, dt.
\end{align*}

Recalling \eqref{starnorm} and \eqref{L2normstarnormgradnorm}, it holds that
\begin{align*}
\|N(\varphi)\|_{L^{2}}  \leq C_{p} \|\nabla N(\varphi)\|_{L^{2}} = C_{p} \|\varphi\|_{*}, \quad \|\varphi\|_{L^{2}}^{2}  \leq \|\varphi\|_{*} \|\nabla \varphi\|_{L^{2}},
\end{align*}
and so we obtain
\begin{subequations}
\begin{align}
\label{singular:ctsdep:K1} |K_{1}| & \leq  \left ( (\mathcal{Y} \lambda_{p} \overline{C} (C_{p}+1)^{2})^{2} + C_{p}^{2} (\mathcal{Y} \lambda_{a} \mathrm{L}_{h} C_{p} )^{2}  \right ) \|\varphi\|_{L^{2}(0,s;*)}^{2} \\
\notag & \quad + \frac{1}{2} \|\nabla \varphi\|_{L^{2}(0,s;L^{2})}^{2}, \\
\label{singular:ctsdep:K2} |K_{2}| & \leq \left ( C_{p}^{2}(\mathcal{Y} \lambda_{p} h_{\infty})^{2} + \left ( \mathcal{Y} \frac{\gamma}{\varepsilon} \|\Lambda''\|_{L^{\infty}(\mathbb{R})} + \frac{\mathcal{Y}^{2}}{2} \right )^{2} \right ) \|\varphi\|_{L^{2}(0,s;*)}^{2}  \\
\notag & \quad + \frac{2 \chi^{2} + 1}{4} \|\sigma\|_{L^{2}(0,s;L^{2})}^{2} + \frac{1}{4} \|\nabla \varphi\|_{L^{2}(0,s;L^{2})}^{2} \\
\notag & \leq \left ( C_{p}^{2}(\mathcal{Y} \lambda_{p} h_{\infty})^{2} + \left ( \mathcal{Y} \frac{\gamma}{\varepsilon} \|\Lambda''\|_{L^{\infty}(\mathbb{R})} + \frac{\mathcal{Y}^{2}}{2} \right )^{2} \right ) \|\varphi\|_{L^{2}(0,s;*)}^{2} \\
\notag & \quad + \frac{(2C_{p}^{2}+1)(2\chi^{2}+1)}{4} \|\nabla \sigma\|_{L^{2}(0,s;L^{2})}^{2} + \frac{1}{4} \|\nabla \varphi\|_{L^{2}(0,s;L^{2})}^{2} + C \|\sigma_{\infty}\|_{L^{2}(0,T;H^{1})}^{2}.
\end{align}
\end{subequations}

We point out that the reason not to use the Poincar\'{e} inequality directly to estimate terms involving $\|\varphi\|_{L^{2}(0,s;L^{2})}^{2}$ in $K_{2}$ is to obtain terms involving $\|\nabla \varphi\|_{L^{2}(0,s;L^{2})}^{2}$ with coefficients that are independent of $\mathcal{Y}$.  Upon collecting all of the terms from \eqref{singular:ctsdep:K3}, \eqref{singular:ctsdep:K4}, \eqref{singular:ctsdep:K1} and \eqref{singular:ctsdep:K2}, and also keeping in mind the estimate \eqref{ctsdep:J1} for $J_{1}$, we find that
\begin{align*}
& \frac{\mathcal{Y}}{2} \|\varphi(s)\|_{*}^{2}  + \mathcal{Z} \frac{\kappa}{4} \|\sigma(s)\|_{L^{2}}^{2} + \left ( \mathcal{Z} \frac{D_{0}}{2} - (2C_{p}^{2}+1) \frac{2 \chi^{2} + 3}{4} - \frac{C_{p}^{2}}{2} \right )\|\nabla \sigma\|_{L^{2}(0,s;L^{2})}^{2} \\
& \qquad + \left ( \mathcal{Y} \gamma \varepsilon - 2 \mathcal{Z} D_{0} \eta^{2} - 2 \overline{C}^{2} \lambda_{c}^{2} \mathcal{Z}^{2}(C_{p}+1)^{2} - \frac{3}{4} \right ) \|\nabla \varphi\|_{L^{2}(0,s;L^{2})}^{2}  \\
& \quad \leq \frac{\mathcal{Y}}{2} \|\varphi_{0}\|_{*}^{2} + \mathcal{Z} \frac{3\kappa}{4} \|\sigma_{0}\|_{L^{2}}^{2} +  \kappa^{2} \mathcal{Z}^{2} \|\partial_{t}\sigma_{\infty}\|_{L^{2}(0,T;L^{2})}^{2} + 2 \kappa \mathcal{Z} \|\sigma_{\infty}\|_{L^{\infty}(0,T;L^{2})}^{2}  \\
& \qquad + C \left ( \mathcal{Y}, \lambda_{p}, \overline{C}, C_{p}, \lambda_{a}, \mathrm{L}_{h}, h_{\infty}, \gamma, \varepsilon, \|\Lambda''\|_{L^{\infty}(\mathbb{R})} \right ) \|\varphi\|_{L^{2}(0,s;*)}^{2} + C \|\sigma_{\infty}\|_{L^{2}(0,T;H^{1})}^{2},
\end{align*}
for some positive constant $C$ not depending on $\kappa$.  Choosing $\mathcal{Z}$ and $\mathcal{Y}$ so that the coefficients on the left-hand side are positive, we obtain
\begin{align*}
& \|\varphi(s)\|_{*}^{2} + \kappa \|\sigma(s)\|_{L^{2}}^{2} + \int_{0}^{s} d_{1} \left ( \|\nabla \sigma\|_{L^{2}}^{2} + \|\nabla \varphi\|_{L^{2}}^{2} \right ) \, dt \\
& \quad \leq d_{2} \left ( \| \varphi_{0}\|_{*}^{2} + \kappa \| \sigma_{0} \|_{L^{2}}^{2} \right ) + d_{3} \|\varphi\|_{L^{2}(0,s;*)}^{2} \\
& \qquad + d_{4} \left ( \|\sigma_{\infty}\|_{L^{2}(0,T;H^{1})}^{2} + \kappa \|\sigma_{\infty}\|_{L^{\infty}(0,T;L^{2})}^{2} + \kappa^{2} \| \sigma_{\infty} \|_{H^{1}(0,T;L^{2})}^{2} \right ),
\end{align*}
for positive constants $d_{1}$, $d_{2}$, $d_{3}$ and $d_{4}$ not depending on $\kappa$.  Application of the Gronwall inequality \eqref{Gronwall} yields the desired result.

\subsection{Full uniqueness}
If $\varphi_{1,0} = \varphi_{2,0}$, $\sigma_{1,0} = \sigma_{2,0}$ and $\sigma_{\infty,1} = \sigma_{\infty,2}$, then it holds that
\begin{align*}
\|\sigma_{1} - \sigma_{2}\|_{L^{2}(\Omega \times (0,T))}^{2} + \|\varphi_{1} - \varphi_{2}\|_{L^{2}(\Omega \times (0,T))}^{2} = 0.
\end{align*}
Hence, $\varphi_{1} = \varphi_{2}$ and $\sigma_{1} = \sigma_{2}$ a.e. in $\Omega \times (0,T)$.  As $\mu_{\infty,1} = \mu_{\infty,2}$, we see that $\mu = \mu_{1} - \mu_{2} \in H^{1}_{0}$. Thus, substituting $\zeta = \mu$ in \eqref{singular:truncated:uniq:varphi}, integrating in time and applying the Poincar\'{e} inequality leads to
\begin{align*}
0 = \|\nabla \mu\|_{L^{2}(0,T;L^{2})}^{2} \geq \frac{1}{C_{p}^{2}} \|\mu\|_{L^{2}(0,T;L^{2})}^{2} \Rightarrow \mu_{1} = \mu_{2} \text{ a.e. in } \Omega \times (0,T).
\end{align*}
Upon setting $\varphi = \mu = \sigma = 0$ in \eqref{singular:truncated:uniq:mu} we obtain
\begin{align*}
\int_{\Omega} \psi \lambda \, dx = 0 \quad \forall \lambda \in H^{1}_{0},
\end{align*}
which in turn implies that $\psi = \psi_{1} - \psi_{2} = 0$ a.e. in $\Omega \times (0,T)$.  This completes the proof of Theorem \ref{thm:singular:ctsdep}.

\section*{Acknowledgments}
The authors would like to thank the referees for their careful reading and their precious suggestions which have improved the manuscript.

\end{document}